\newcommand{\xBc}{\langle}
\newcommand{\xBe}{\rangle}

\newcommand{\xbD}{\Delta}
\newcommand{\xbF}{\Phi}
\newcommand{\xbG}{\Gamma}

\newcommand{\xba}{\alpha}
\newcommand{\xbb}{\beta}

\newcommand{\xbe}{\in}
\newcommand{\xbf}{\phi}

\newcommand{\xbm}{\mu}

\newcommand{\xbo}{\omega}
\newcommand{\xbp}{\pi}

\newcommand{\xbs}{\sigma}
\newcommand{\xbt}{\tau}
\newcommand{\xbx}{\xi}

\newcommand{\xCN}{\neg}
\newcommand{\xCO}{ }

\newcommand{\xCf}{\hspace{0.1em}}

\newcommand{\xcA}{\forall}

\newcommand{\xcE}{\exists}

\newcommand{\xcH}{\not\Rightarrow}

\newcommand{\xcO}{\bigvee}
\newcommand{\xcP}{\not\rightarrow}

\newcommand{\xcU}{\bigwedge}

\newcommand{\xcc}{\subseteq}

\newcommand{\xce}{\not\in}

\newcommand{\xcg}{\geq}
\newcommand{\xch}{\Rightarrow}

\newcommand{\xcj}{\Leftrightarrow}
\newcommand{\xck}{\leq}

\newcommand{\xcm}{\models}

\newcommand{\xco}{\vee}
\newcommand{\xcp}{\rightarrow}

\newcommand{\xcr}{\leftrightarrow}

\newcommand{\xcu}{\wedge}
\newcommand{\xcv}{\cup}

\newcommand{\xcz}{\Box}

\newcommand{\xDH}{\item }

\newcommand{\xDO}{\circ}

\newcommand{\xdF}{\mbox{\boldmath$F$}}

\newcommand{\xdT}{\mbox{\boldmath$T$}}

\newcommand{\xdu}{{\cal U}}

\newcommand{\xEI}{\begin{itemize}}
\newcommand{\xEJ}{\end{itemize}}

\newcommand{\xET}{\%}

\newcommand{\xEd}{\neq}

\newcommand{\xEh}{\begin{enumerate}}

\newcommand{\xEj}{\end{enumerate}}

\newcommand{\xEn}{\begin{description}}
\newcommand{\xEp}{\end{description}}

\newcommand{\xeH}{\stackrel{\rightarrow}{\nrightarrow}}

\newcommand{\xeb}{\prec}

\newcommand{\xex}{\upharpoonright}

\newcommand{\xfA}{\mid}

\newcommand{\Xl}{\ldots}

\newcommand{\bl}{\begin{lemma} \rm}
\newcommand{\el}{\end{lemma}}
\newcommand{\br}{\begin{remark} \rm}
\newcommand{\er}{\end{remark}}
\newcommand{\be}{\begin{example} \rm}
\newcommand{\ee}{\end{example}}
\newcommand{\bco}{\begin{corollary} \rm}
\newcommand{\eco}{\end{corollary}}
\newcommand{\bc}{\begin{claim} \rm}
\newcommand{\ec}{\end{claim}}
\newcommand{\bfa}{\begin{fact} \rm}
\newcommand{\efa}{\end{fact}}
\newcommand{\bp}{\begin{proposition} \rm}
\newcommand{\ep}{\end{proposition}}
\newcommand{\bd}{\begin{definition} \rm}
\newcommand{\ed}{\end{definition}}
\newcommand{\bcs}{\begin{construction} \rm}
\newcommand{\ecs}{\end{construction}}
\newcommand{\bcd}{\begin{condition} \rm}
\newcommand{\ecd}{\end{condition}}
\newcommand{\bt}{\begin{theorem} \rm}
\newcommand{\et}{\end{theorem}}
\newcommand{\bn}{\begin{notation} \rm}
\newcommand{\en}{\end{notation}}
\newcommand{\bfi}{\begin{bild} \rm}
\newcommand{\efi}{\end{bild}}
\newcommand{\bsta}{\begin{statement} \rm}
\newcommand{\esta}{\end{statement}}
\newcommand{\bcom}{\begin{comment} \rm}
\newcommand{\ecom}{\end{comment}}
\newcommand{\bdia}{\begin{diagram} \rm}
\newcommand{\edia}{\end{diagram}}

\newcommand{\bfc}{\begin{figure}[htb] \begin{center}}
\newcommand{\efc}{\end{center} \end{figure}}

\documentclass{article}

\usepackage{amssymb,latexsym,epic,eepic,rotating}

\title{Further Comments on Yablo's Construction
\thanks{File: WAB, [Sch25], arXiv 2504.10370
}
}

\author{Karl Schlechta
\thanks{
schcsg@gmail.com - https://sites.google.com/site/schlechtakarl/ -
Koppeweg 24, D-97833 Frammersbach, Germany}
\thanks{
Retired, formerly: Aix-Marseille Universit\'{e}, CNRS, LIF UMR 7279, F-13000
Marseille, France
}
}

\begin{document}

\newtheorem{lemma}{Lemma}[section]
\newtheorem{theorem}[lemma]{Theorem}
\newtheorem{proposition}[lemma]{Proposition}
\newtheorem{corollary}[lemma]{Corollary}
\newtheorem{claim}[lemma]{Claim}
\newtheorem{fact}[lemma]{Fact}
\newtheorem{remark}[lemma]{Remark}
\newtheorem{definition}{Definition}[section]
\newtheorem{construction}{Construction}[section]
\newtheorem{condition}{Condition}[section]
\newtheorem{example}{Example}[section]
\newtheorem{notation}{Notation}[section]
\newtheorem{bild}{Figure}[section]
\newtheorem{comment}{Comment}[section]
\newtheorem{statement}{Statement}[section]
\newtheorem{diagram}{Diagram}[section]

\renewcommand{\labelenumi}
  {(\arabic{enumi})}
\renewcommand{\labelenumii}
  {(\arabic{enumi}.\arabic{enumii})}
\renewcommand{\labelenumiii}
  {(\arabic{enumi}.\arabic{enumii}.\arabic{enumiii})}
\renewcommand{\labelenumiv}
  {(\arabic{enumi}.\arabic{enumii}.\arabic{enumiii}.\arabic{enumiv})}

\maketitle

\setcounter{secnumdepth}{3}
\setcounter{tocdepth}{3}

\begin{abstract}

The paper shows that any paradoxical construction of the same general
kind as Yablo's actually contains his construction within it,
in a sense that can be made mathematically precise.
More specifically, if we represent Yablo's construction and another
of the same general kind by labelled acyclic graphs satisfying
certain conditions, then the Yablo graph is always a subgraph of the
other one.

We use this result for a characterisation of this kind of labelled
acyclic graphss.

\end{abstract}

\tableofcontents
\clearpage

% *** BEGIN LATEX SOURCE wab-1.tex ***
%
% Aus Karltex File: wab-1.ms
%
%

$ \xCO $
% {\tiny KS= wab-Anf   FOUND IN START
% c:-ks-w-t-framart-phil-yablo-wablo-wba-1.mp} \\[1mm]
\clearpage
\section{
Introduction
}

% {\tiny LABEL: {Section Introduction}} \\[1mm]
\label{Section Introduction}

Yablo's construction is, essentially, the result of ``uncoiling''
the liars' paradox into an infinite sequence without end of
contradictions. When the nodes in this construction represent
propositional variables, any truth value attributed to some
fixed node leads to a contradiction ``down the line''.

The question is to characterise all directed acyclic graphs with
this property about truth values (and some additional properties).
We do this by showing that any such graph contains Yablo's
construction as a subgraph.

The basic idea is to show that such graphs must contain an interplay
of local contradictions and paths to the next contradiction,
properties (C1) and (C2) below, see
Section \ref{Section C1C2} (page \pageref{Section C1C2}).

More precisely, we first show that basic contradictions all (roughly) have
the form of triangles, $x,$ $y,$ $z,$ where $x$ says $ \xCN y$ and $ \xCN
z,$ and $y$ says $ \xCN z,$ thus,
if we go directly from $x$ to $z,$ then $x$ says $ \xCN z,$ but if we take
the detour
via $y,$ then $x$ says $z.$ This is how Yablo's structure works.
We then show that there is one essential way to put those triangles
together,
connecting them via $y$ (in above example). This is how Yablo's original
structure works. We may in addition do something
via $z,$ but this alone will not be sufficient.

Thus, any such structure contains Yablo's structure, whether we construct
it from scratch, or find it within an existing graph.
\subsection{
Outline
}

% {\tiny LABEL: {Section Outline}} \\[1mm]
\label{Section Outline}

These notes are based on  \cite{Sch22} and
 \cite{Sch23b}, the reader may have copies of both
ready for more examples and discussions.
See also  \cite{Sch25}.
(Please note that arXiv pdf files sometimes have problems with some of
the diagrams. In that case, please consult the html version.)
We liberally use material from these sources.
All these papers are, of course, based on Yablo's seminal
 \cite{Yab82}.

Yablo gave an (infinite, propositional) theory using variables $x_{i},$
$i< \xbo $
such that neither $x_{0}$ nor $ \xCN x_{0}$
has a model. He showed how to destroy all attempted models. A successful
model $m$ begins with $m \xcm x_{0}$ (or $m \xcm \xCN x_{0}),$ and
continues in a coherent way for
all $x_{i},$ $0<i< \xbo.$ For each such $m,$ there is in Yablo's
construction some $i_{m}$
such that $m$ cannot be continued beyond $i_{m}$ - there is a
contradiction
starting at $i_{m}.$

We look at two questions:
 \xEh
 \xDH Are there contradictions other than Yablo's?
 \xDH are there other ways to organize the destructions than those chosen
by Yablo?
 \xEj

The answers are:
 \xEh
 \xDH Essentially not,
see Section 
\ref{Section Elementary-Cells} (page 
\pageref{Section Elementary-Cells}),
 \xDH yes, but they contain Yablo's organisation,
see Section 
\ref{Section Combining-Cells} (page 
\pageref{Section Combining-Cells}), and
Section 
\ref{Section Finite-Infinite} (page 
\pageref{Section Finite-Infinite})  and
Section \ref{Section C1C2} (page \pageref{Section C1C2})  for an overview.
 \xEj

Thus, Yablo's solution is essentially universal.

Prerequisites:

 \xEh
 \xDH
All graphs considered will by acyclic.

 \xDH
In most of the paper, formulas will be of the type $ \xcU \xCN x_{i},$ in
Section \ref{Section Result} (page \pageref{Section Result})
we also consider trivial formulas of the type $(y \xco \xCN y \xco \xdT
),$
see Remark \ref{Remark Neutral} (page \pageref{Remark Neutral}).

 \xDH For ease of exposition, we introduce a third truth value, $ \xbx,$
as abbreviation for ``$ \xbf $ and $ \xCN \xbf $ are inconsistent''.
See Definition \ref{Definition Graph} (page \pageref{Definition Graph}), (4.4).

 \xDH We will use only one type of contradictory cell
(see Definition 
\ref{Definition Yablo-Structure} (page 
\pageref{Definition Yablo-Structure}), 2.1)
in a graph.
This prevents ``cheating'' (where the real work is done by some
more complicated structure)
and focusses attention on the type considered, see
Definition 
\ref{Definition Cell-Hierarchy} (page 
\pageref{Definition Cell-Hierarchy}), (1),
and Remark \ref{Remark ESC-Simple} (page \pageref{Remark ESC-Simple}).

 \xDH Contradictions in a graph will begin immediately, and not after some
arbitrary path, except for
Section \ref{Section PathCells} (page \pageref{Section PathCells}).
This simplifies constructions, and is not an important
restriction.

 \xDH We do not admit any additional theory outside the graph.
(See Section 
\ref{Section Diamond-Problem} (page 
\pageref{Section Diamond-Problem}).)
 \xEj
\subsubsection{
Overview of the Paper
}

% {\tiny LABEL: {Section Overview}} \\[1mm]
\label{Section Overview}

 \xEh

 \xDH Section 
\ref{Section Logic} (page 
\pageref{Section Logic})  is sowhat exceptional,
we discuss the relation of the (intuitive) logic
we use here, to classical propositional logic. We also outline briefly
a formal logic adapted to path-based reasoning in directed graphs.

This section is a sketch only, its purpose is to indicate how to
make allusions to a background logic more formal. It is not central to
the main ideas and results of the paper.

After this section, we will just say ``classical logic'', and mean
a suitably weakend version, as indicated here.

 \xDH In Section \ref{Section Basics} (page \pageref{Section Basics}), we
first present the background, its logical and graph theoretical side,
and their interplay. We also discuss a useful third truth value, which
expresses that neither $ \xbf $ nor $ \xCN \xbf $ is consistent.

 \xDH In Section \ref{Section Strategy} (page \pageref{Section Strategy}),
we discuss our strategy, and the
antagonistic conditions (C1) and (C2)
and how to satisfy them in an inductive construction.
This section is a preparation to
Section \ref{Section Combining-Cells} (page \pageref{Section Combining-Cells}).

 \xDH In Section 
\ref{Section Comparison} (page 
\pageref{Section Comparison}), we compare our results to
those
of related papers.

 \xDH In Section 
\ref{Section Elementary-Cells} (page 
\pageref{Section Elementary-Cells}),
we discuss elementary contradictions (cells). Essentially, these are
variations of Yablo's triangles. The trivial contradiction, see
Definition 
\ref{Definition Cell-Hierarchy} (page 
\pageref{Definition Cell-Hierarchy}), is too
simple, it has an ``escape possibility'', i.e., combining them will
leave one possibility without contradictions. See
Definition \ref{Definition Escape} (page \pageref{Definition Escape}).
Diamonds, constructions with
four sides, will not work because they need a ``synchronization'', which
is not available in our framework.
See Section 
\ref{Section Diamond-Problem} (page 
\pageref{Section Diamond-Problem}).
We also show that all sides of the triangles have to be negative, see
Remark \ref{Remark Simple-Cells} (page \pageref{Remark Simple-Cells}).
Section \ref{Section YC-Details} (page \pageref{Section YC-Details})
discusses further details of these triangles.
We generalize to triangles formed by paths in
Section \ref{Section PathCells} (page \pageref{Section PathCells}).

 \xDH Section 
\ref{Section Combining-Cells} (page 
\pageref{Section Combining-Cells})
is the central part of the paper,
Section 
\ref{Section Result} (page 
\pageref{Section Result})  will essentially be a corollary.
Here, we discuss how to put contradiction cells together to obtain graphs
which code the liar paradox.

We discuss 4 potential possibilities, show that one does not work,
but the other three all are or contain necessarily Yablo's construction.

The constructions are necessary to obtain Yablo's result, whether
we construct a graph from scratch, or (try to) construct a suitable
subgraph within a given graph.
This is the main argument and result of the paper, the final result,
Section 
\ref{Section Result} (page 
\pageref{Section Result})  is then an easy consequence.

More precisely, we first follow Yablo's construction, see
Section \ref{Section Combining-Yablo} (page \pageref{Section Combining-Yablo}).

We then generalize slightly,
see
Section 
\ref{Section Combining-Saw-Blade} (page 
\pageref{Section Combining-Saw-Blade})
and the discussion of ``Saw Blades'' in
Section 7.6 of  \cite{Sch22} and
Section 3 in  \cite{Sch23b},
but show that this construction contains Yablo's
construction, so the latter is minimal.

In both cases, we show that the construction is necessary.

 \xDH
In Section \ref{Section Result} (page \pageref{Section Result}),
we show that
 \xEh
 \xDH a graph without valuation $ \xbG $ can be given a valuation such
that there
is an element $x$ such that neither $ \xbf_{x}$ nor $ \xCN \xbf_{x}$ are
consistent, if
there is a suitable injection from Yablo's structure into $ \xbG,$ and
a partial valuation for the nodes of the paths of the injection, such
that contradictions are preserved. We complete the partial injection
by a trivial extension.
 \xDH conversely, if there is a graph with (complete) valuation with some
node $x$ such that neither $ \xbf_{x}$ nor $ \xCN \xbf_{x}$ are
consistent, then we can
construct a subgraph which corresponds again to an injection of the
Yablo structure. We use here again the minimality of Yablo's construction,
this time imside the given graph.
 \xEj
The minimality of Yablo's construction is the core of the argument
in both directions.
 \xEj
\clearpage
\subsection{
A Comment on Logic
}

% {\tiny LABEL: {Section Logic}} \\[1mm]
\label{Section Logic}

We take here a closer look at the background logic of Yablo's and
similar constructions. The reader will not miss much if he trusts his
intuition, just think: ``a somewhat weakened, possibly infinite,
classical logic, following the arrows in the graph''.

As this section is somewhat informal, we refer the reader to
Section \ref{Section Basics} (page \pageref{Section Basics})  for definitions.

For a short comparison of the present section to the work of
 \cite{RRM13}, see
Section \ref{Section Logics-Rabern} (page \pageref{Section Logics-Rabern}).

The present section has conceptually two parts. The first part,
Section 
\ref{Section Yablo-Requirements} (page 
\pageref{Section Yablo-Requirements}),
discusses
how to find inconsistencies, but this is perhaps a too simple minded
logic.

The second part,
Section \ref{Section Rules} (page \pageref{Section Rules}),
discusses more general questions, e.g., also how to
treat partial inconsistencies. We present here basically only a problem,
which needs attention, and leave the subject to another paper.
\subsubsection{
Fundamentals
}

Intuitively speaking,
Yablo's example is formulated in classical propositional logic, but
admitting
infinite formulas. (Infinite formulas are required.) When we speak about
finite fragments of the graph, the logic is simply classical. All this
with some caveats, classical logic has more consequences than we want.

We make this more precise now.

We have two possibilities,
we may approach such a logic from above (e.g. classical logic), or from
below.
Building from below may be easier, but may be too
specific to merit the label ``logic''.
We sketch here some ideas for the approach from below.

We assume a (directed - as in Yablo's case - or not) graph, with formulas
(in a
given language) attached to the nodes, as well as formulas describing a
logical
connection between the formulas of the nodes and their (graph) neighbours.
In
Yablo's case, the formulas attached to the nodes are propositional
variables,
and the logical connections are of the type $x \xcr \xcU \xCN x_{i}.$

\bd

$\hspace{0.01em}$

% {\tiny (+++ Orig. No.:  Definition Fundamental Axiom +++)}

% {\tiny LABEL: {Definition Fundamental Axiom}} \\[1mm]
\label{Definition Fundamental Axiom}

(Fundamental Axiom)

Relevance is established ONLY via paths in the graph.

Thus, there is nothing like EFQ.

If there is no path from $ \xCf a$ to $b,$ then there is no connection
between the two. Also, e.g., if the graph consists of
disconneted parts, then variables in one part have nothing to do
with variables in another part, even if they have the same name.
\subsubsection{
Requirements for Yablo-like structures
}

% {\tiny LABEL: {Section Yablo-Requirements}} \\[1mm]
\label{Section Yablo-Requirements}

\ed

What we really need for Yablo's argument is forward reasoning along the
paths defined by the graph, and backtracking when we found a contradiction
by forward reasoning, i.e. found $y$ and $ \xCN y$ following some paths
originating
in $x.$ In that case we will conclude that $x$ is impossible, likewise
for $ \xCN x$ when we started with $ \xCN x.$
The forward part is relatively straightforward, essentially, we ``collect
''
the consequences along the arrows,
respecting the classical laws about $ \xcu,$ $ \xco,$ $ \xCN.$
We will illustrate this in
Example \ref{Example Logic} (page \pageref{Example Logic})  below,
generalized in (7) there.

General backtracking is more complicated, see e.g.
Example \ref{Example Origin} (page \pageref{Example Origin}).

We must not forget to add that there are no other rules, in particular
we may not close classically.

Note also, that, if we conclude e.g. that the assumption $ \xfA x \xfA =
\xdT $ leads
to a contradiction, we must not necessarily conlude $ \xfA x \xfA = \xdF
,$ but
only that $ \xfA x \xfA = \xdT $ is impossible, both $ \xfA x \xfA = \xdT
$ and $ \xfA x \xfA = \xdF $ may
be impossible - as in Yablo's structure.

Finally we define

\bd

$\hspace{0.01em}$

% {\tiny (+++ Orig. No.:  Definition Model +++)}

% {\tiny LABEL: {Definition Model}} \\[1mm]
\label{Definition Model}

A model is a coherent assignment of classical truth values according to
the
laws of the truth values and the operators.

E.g., for $x_{1} \xcr \xCN x_{2} \xcu \xCN x_{3},$ $x_{2} \xcr \xCN
x_{3},$
$ \xfA x_{1} \xfA = \xfA x_{2} \xfA = \xdF,$ $ \xfA x_{3} \xfA = \xdT $
is a model, but there is no model with $ \xfA x_{1} \xfA = \xdT.$

$ \xfA x_{1} \xfA = \xfA x_{2} \xfA = \xfA x_{3} \xfA = \xbx $ is not a
model, though it is coherent.

\ed

\be

$\hspace{0.01em}$

% {\tiny (+++ Orig. No.:  Example Logic +++)}

% {\tiny LABEL: {Example Logic}} \\[1mm]
\label{Example Logic}

We now discuss some examples, where it seems evident what to do, and
describe models, where possible.
Thus, we have some examples any such logic must satisfy.

 \xEh

 \xDH $x \xcP y \xcP z,$ $x \xcP z,$ $x= \xCN y \xcu \xCN z,$ $y= \xCN z$

 \xEh
 \xDH $ \xfA x \xfA = \xdT:$ This assumption leads to a contradiction,
there are no models.

 \xDH $ \xfA x \xfA = \xdF:$ Thus, $ \xCN x=y \xco z,$ $y= \xCN z.$

We have two models:

$m_{1}:$ $ \xCN x,$ $y,$ $ \xCN z,$ $m_{2}:$ $ \xCN x,$ $ \xCN y,$ $z$

 \xEj

 \xDH $x \xcP y \xcP z,$ $x \xcP z,$ $x \xcP y',$ $x= \xCN y \xcu \xCN z
\xcu \xCN y',$ $y= \xCN z$

 \xEh

 \xDH $ \xfA x \xfA = \xdT:$ This assumption leads to a contradiction,
there are no models.

Note: If we have a contradiction in part of a
(finite or infinite, see (5) below)
conjunction, then the whole conjunction is contradicted.

 \xDH $ \xfA x \xfA = \xdF:$ Thus, $ \xCN x=y \xco z \xco y'.$

We have 4 models: $m_{1}:$ $ \xCN x,$ $y,$ $ \xCN z,$ $y',$ $m_{2}:$ $
\xCN x,$ $y,$ $ \xCN z,$ $ \xCN y',$
$m_{3}:$ $ \xCN x,$ $ \xCN y,$ $z,$ $y',$ $m_{4}:$ $ \xCN x,$ $ \xCN y,$
$z,$ $ \xCN y'.$

 \xEj

 \xDH $x \xcP y \xcP z,$ $x \xcP z,$ $y \xcP z' $ $x= \xCN y \xcu \xCN z,$
$y= \xCN z \xcu \xCN z' $

We open a new possibility for $ \xCN y.$

 \xEh

 \xDH $ \xfA x \xfA = \xdT:$ We have 1 model: $m:$ $x,$ $ \xCN y,$ $ \xCN
z,$ $z' $

Note: if we have a contradiction in ALL parts of a
(finite or infinite, see (5) below)
disjunction, then
the disjunction is contradicted, but not if a part of the disjunction
is not contradicted.

 \xDH $ \xfA x \xfA = \xdF:$ We have 3 models:

$m_{1}:$ $ \xCN x,$ $y,$ $ \xCN z,$ $ \xCN z',$ $m_{2}:$ $ \xCN x,$ $
\xCN y,$ $z,$ $z',$ $m_{3}:$ $ \xCN x,$ $ \xCN y,$ $z,$ $ \xCN z'.$

 \xEj

 \xDH $x \xcP y \xcP z,$ $x \xcP z,$ $x \xcP z',$ $y \xcP z' $ $x= \xCN y
\xcu \xCN z \xcu \xCN z',$ $y= \xCN z \xcu \xCN z' $

 \xEh

 \xDH $ \xfA x \xfA = \xdT:$ This assumption leads to a contradiction,
there are no models.

We closed the possibility opened in (3).

 \xDH $ \xfA x \xfA = \xdF:$ We have 4 models:

$m_{1}:$ $ \xCN x,$ $y$ $, \xCN z,$ $ \xCN z',$ $m_{2}:$ $ \xCN x,$ $
\xCN y$ $,z,$ $z',$ $m_{3}:$ $ \xCN x,$ $ \xCN y$, $z,$ $ \xCN z',$
$m_{4}:$ $ \xCN x,$ $ \xCN y$, $ \xCN z,$ $z'.$

 \xEj

 \xDH Yablo's Construction

 \xEh
 \xDH $ \xfA x \xfA = \xdT:$ Take any $x_{i}.$ Take any $x_{i} \xcP
x_{j}$ (recall that one such $x_{j}$ would
suffice), take any $x_{j} \xcP x_{k}$ (recall that here, all $x_{k}$ are
needed), then
there is $x_{i} \xcP x_{k},$ and we are done.

 \xDH $ \xfA x \xfA = \xdF:$ It suffices to show the result for $x_{0},$
but the same result
and argument holds for any $x_{i},$ so we start at $x_{i}.$ As we have an
``OR'' at
$x_{i}$ for $ \xfA x_{i} \xfA = \xdF,$ for any $x_{i} \xcP x_{j},$ $ \xfA
x_{j} \xfA = \xdT $ must be impossible, as was
shown above.

 \xEj

 \xDH Finite fragment of Yablo's construction

Consider the Yablo construction up to $x_{k},$ i.e. $\{x_{i}:i \xck k\},$
then $ \xfA x_{k} \xfA = \xdT $
fails to be contradictory, and thus $ \xfA x_{0} \xfA = \xdF $ fails to be
contradictory, too.
(As a matter of fact, we may fill in truth values from $x_{k}$ downward.

 \xDH We generalize:

In the following, $ \xcA $ has to be replaced by $ \xcU,$ etc.

 \xEh

 \xDH We just formalize the rules we have found, the possible
contradictions,
essentially the rules (C1) and (C2).

 \xDH
If $ \xcE (x_{i} \xcP x_{j})$ ($ \xcA (x_{j} \xcP x_{k})$ $ \xcE (x_{i}
\xcP x_{k})$), then $x_{i}$ cannot be true

 \xDH
If the Yablo triangles are formed by paths instead of simple arrows,
we have to generalize suitably.
But the principle is the same, we write the construction down.

 \xDH If
$ \xcA (x_{0} \xcP x_{i})$ $x_{i}= \xdF,$ then $ \xCN x_{0}$ cannot be
true, i.e.

if
$ \xcA (x_{0} \xcP x_{i})$ ($ \xcE (x_{i} \xcP x_{j})$ ($ \xcA (x_{j}
\xcP x_{k})$ $ \xcE (x_{i} \xcP x_{k})$)), then $ \xCN x_{0}$ cannot be
true.

 \xEj

 \xEj
\subsubsection{
A more ambitious project
}

% {\tiny LABEL: {Section Rules}} \\[1mm]
\label{Section Rules}

\ee

 \xEh

 \xDH Framework

We try to simplify the sets of formulas corresponding to the given
graph, to be able to treat partial results.
(For sake of simplicity, we omit the truth value $ \xbx $
and work with the definition.)

 \xDH General rules

 \xEh
 \xDH Consider the formula $ \xdT \xco \xbf \xco \xCN \xbf $ used for
clipping unwanted
details (see Section \ref{Section Result} (page \pageref{Section Result})).

Logically, this is, of course, equivalent to $ \xdT $ (even if $ \xfA \xbf
\xfA = \xbx),$ but
we added the $ \xbf $-part to make the formula correspond to the graph.

So, we may replace $ \xdT \xco \xbf \xco \xCN \xbf $ by $ \xdT.$ Note
that this is a logical
equivalence, but we will not replace $ \xdT $ by $ \xdT \xco \xbf \xco
\xCN \xbf,$ as this would
violate the correspondence between formulas and graph in an
essential way (not forgetting the correspondence, but inventing
something).

We write $ \xdT \xeb \xdT \xco \xbf \xco \xCN \xbf $ to indicate the
formulas should
only be read as going from complicated to simpler.

 \xDH Of course, it is more elegant to use a more general principle: $
\xdT \xeb \xdT \xco \xbf.$

Likewise, we may introduce the principle $ \xbf \xeb \xbf \xcu \xdT.$

These principles seem basic.

 \xDH A more complicated case is the following (Yablo cell):

$x \xEd \xdT $ $ \xeb $ $((x \xcr \xCN y \xcu \xCN z)$ $ \xcu $ $(y \xcr
\xCN z)).$

 \xEj

If we want to treat partial contradictions, we
have to do some careful
bookkeeping, as we argue with arrows, and not with nodes. In particular,
we must take care not to use the same node too many times, as the
following
eaxmple shows:

\be

$\hspace{0.01em}$

% {\tiny (+++ Orig. No.:  Example Origin +++)}

% {\tiny LABEL: {Example Origin}} \\[1mm]
\label{Example Origin}

(It shows that we cannot neglect the point of origin.)

Let $x \xcP y \xcP z,$ $x' \xcP y \xcP z',$ $x \xcP z,$ $x' \xcP z'.$
Add now $x' \xcP z.$

We have two cells joined at the knee, one contradiction, with head $x',$
is complete, the other is not. So for the complete contradiction, we may
forget $y \xcP z$ and $y \xcP z',$ and set $ \xfA x' \xfA \xEd \xdT,$
but not for the incomplete one.

So, we have to memorize the node relative to which $y \xcP z$ and $y \xcP
z' $
are contradicted. Here, $y \xcP z' $ is contradicted relative to $x',$
but not
relative to $x.$ $y \xcP z$ is contradicted relative to both $x$ and $x'
.$

Thus, the structure plays an essential role, unsurprisingly.

\ee

For illustration, try to follow
Construction 
\ref{Construction Details} (page 
\pageref{Construction Details})  and
Diagram \ref{Diagram Induction} (page \pageref{Diagram Induction})  in
Section \ref{Section Combining-Cells} (page \pageref{Section Combining-Cells}).

 \xEj
\subsubsection{
Summary and further comments
}

From now on, when we say classical logic, we mean a suitably impoverished
version of classical logic, or a suitably developed logic for the
approach from below.
In this sense, the logical problems in this paper are trivial. The real
problems
are combinatorial, analysing the essential patters in the graphs.

Further comments:

 \xEh

 \xDH
We will see how a graph with valuation defines a (propositional) theory
(see Definition \ref{Definition Graph} (page \pageref{Definition Graph}), (3)).
Thus, when we discuss different such graphs, we have
different theories. An important part of the present article discusses
building Yablo-like graphs inductively
(see Section 
\ref{Section Elementary-Cells} (page 
\pageref{Section Elementary-Cells})  and
Section 
\ref{Section Combining-Cells} (page 
\pageref{Section Combining-Cells})),
so we discuss
different corresponding theories, and their inductive development,
including
problems in the construction (escape possibilities,
see Definition \ref{Definition Escape} (page \pageref{Definition Escape})).
In particular, we do not only discuss the theory defined by Yablo's
example.

 \xDH
We shall introduce a third truth value, $ \xbx,$ which is only
an abbreviation: given a (classical) theory of a graph, a node $x$ has
truth
value $ \xbx $ iff the assumption that $x$ is true leads to a (classical)
contradiction,
and so does the assumption that $x$ is false. Of course, in finite
classical
propositional logic, this is impossible, but in the infinite version it is
possible, as Yablo has shown (in Yablo's example, every node has truth
value $ \xbx).$

We show certain qualities of $ \xbx,$ which are directly
derived from the definition (and classical logic).
As pointed out by a referee, the properties of $ \xbx $
are the same as those of $ \xdu $ in Kleene's strong 3-valued
logic (it is easily shown that this must be so, but for our purposes it is
purely coincidental).
See Remark \ref{Remark Xi} (page \pageref{Remark Xi}).

 \xEj

To make these considerations more precise would be a different paper.
\clearpage
\subsection{
Basic Definitions and Facts
}

% {\tiny LABEL: {Section Basics}} \\[1mm]
\label{Section Basics}

\bd

$\hspace{0.01em}$

% {\tiny (+++ Orig. No.:  Definition Graph +++)}

% {\tiny LABEL: {Definition Graph}} \\[1mm]
\label{Definition Graph}

 \xEh
 \xDH Background logic

The background logic is classical propositional logic, with perhaps
infinite formulas - with the caveat formulated above in
Section \ref{Section Logic} (page \pageref{Section Logic}).

All notions like ``contradiction'', ``consequence'' etc. are
with respect to this logic.

Implications will be written $ \xch,$ $ \xcj,$ etc., their negations $
\xcH,$ etc.

 \xDH Graphs
 \xEh
 \xDH We consider directed, acyclic graphs, often noted $ \xbG $ etc.
They will usually have one root, often called $x_{0}.$
 \xDH Paths will be written $x-y-z$, $ \xCf x \Xl.y \Xl.z$, etc. (read
from left to right).
 \xDH The nodes will have labels, labels will be unique in a given graph.
For simplicity, we identify labels with the (names of the) nodes.
The labels will be the propositional variables in the (propositional)
language defined by the graph.
 \xEj
 \xDH Valuations
 \xEh
 \xDH We often associate a formula $ \xbf_{x}$ with a node $x,$ such that
the variables
contained in $ \xbf_{x}$ are exactly the direct successors of $x$ in the
graph.
The meaning is that $x$ and $ \xbf_{x}$ are logically equivalent, $x \xcj
\xbf_{x}.$

For simplicity, we often write $x= \xbf $ instead of $ \xbf_{x}= \xbf.$
 \xDH The set of formulas $ \xbf_{x},$ $x$ a node in $ \xbG $ will be
called a valuation of
$ \xbG,$ or the theory of $ \xbG,$ $Th(\xbG),$ given such a valuation.
So $Th(\xbG)$ is a set of equivalences, $\{x \xcj \xbf_{x}:$ $x$ a node
in $ \xbG \}.$
Of course, usually $ \xbG $ will have many different valuations.

 \xDH Most of the time, valuations in this article will be of the
syntactic form
$x \xcj \xcU \{ \xCN x_{i}:i \xbe I\} \xcu \xcU \{x_{j}:j \xbe J\}$ or
even just $x \xcj \xcU \{ \xCN x_{i}:i \xbe I\}.$
In this case, we may write the valuation directly into the graph
by $\{x \xcP x_{i}:i \xbe I\} \xcv \{x \xcp x_{j}:i \xbe J\}$ or just $\{x
\xcP x_{i}:i \xbe I\}.$ This increases
readability considerably.

 \xEj

 \xDH Contradictions and truth values

Let $x$ be a node in a fixed graph $ \xbG $ with valuation as above.

 \xEh

 \xDH Additional assumptions:

In addition to $Th(\xbG),$ we assume sometimes that some nodes have
truth
values, to test the behaviour of certain nodes relative to these
assumptions.

 \xDH Truth values:
 \xEh
 \xDH We write $ \xfA x \xfA = \xdT $ if we set $x$ true, or if it turns
out that $x$ is true
(perhaps under some assumptions). Likewise $ \xfA x \xfA = \xdF.$

More formally, $ \xfA x \xfA = \xdT $ iff for $Th(\xbG)$ and perhaps
some assumptions $ \xfA x_{i} \xfA = \xdT $ or
$= \xdF,$ $Th(\xbG) \xcv \{ \xfA x_{i} \xfA = \xdT $ or $= \xdF:$ $i
\xbe I\}$ implies $ \xfA x \xfA = \xdT,$ etc.

Examples:
 \xEh
 \xDH If $ \xbG $ consists of $x \xcP y,$ with the valuation $x \xcj \xCN
y,$ then the
assumption $ \xfA x \xfA = \xdT $ results in $ \xfA y \xfA = \xdF.$

 \xDH If $ \xbG $ consists of $x \xcP y,$ $x \xcP z,$ $y \xcP z,$ with the
valuations
$x \xcj \xCN y \xcu \xCN z,$ $y \xcj \xCN z,$ then both assumptions $ \xfA
z \xfA = \xdT $ and $ \xfA z \xfA = \xdF $ result in
$ \xfA x \xfA = \xdF.$
 \xEj
 \xEj
 \xDH Contradictions

We say that $x$ is contradictory in some $ \xbG $ with fixed valuation
iff, for all
other possible assumptions, the assumption $ \xfA x \xfA = \xdT $ leads to
a contradiction.

In the second example above, $x$ is contradictory.

If $ \xbG $ is Yablo's structure, then for any node $x_{i}$ both $ \xfA
x_{i} \xfA = \xdT $ and $ \xfA x_{i} \xfA = \xdF $
are contradictory.

 \xDH The third truth value, $ \xbx $

We introduce a third truth value $ \xbx,$ which is an abbreviation
for $x$ and $ \xCN x$ are contradictory - relative to $ \xbG $ with
valuation, of course.
See Remark \ref{Remark Xi} (page \pageref{Remark Xi}).

So, in Yablo's structure, all $x_{i}$ have truth value $ \xbx,$ and,
conversely, if
we want to assume that some $x$ has truth value $ \xbx,$ we may think of
setting
$x=x_{0}$ in Yablo's structure.

 \xDH Simplifications of notation

When the context is clear, we sometimes write $x= \xdT $ instead of $ \xfA
x \xfA = \xdT,$ etc.

As the valuations will (mostly) have the form $ \xcU \xCN x_{i},$
in the case $x= \xdT,$ we have a
conjunction, and in the case $x= \xdF,$ we have a disjunction. This is
trivial,
but very important.

 \xEj

 \xDH Contradictions and paths:

(See also Section \ref{Section Logic} (page \pageref{Section Logic}).)

 \xEh
 \xDH
Contradictions via paths have in a graph the following form: there is a
node
$x$ and a node $y,$ and two paths from $x$ to $y,$ which branch at $x$ and
join
again at $y,$ and they contradict each other, i.e. one branch results in
$y,$
the other in $ \xCN y.$

These are the basic contradictions.

Branchings (more precisely branchings in the logical form of ``OR'') open
new possibilities to avoid contradictions.
See Section 
\ref{Section Combining-Cells} (page 
\pageref{Section Combining-Cells}), in
particular Diagram \ref{Diagram Induction} (page \pageref{Diagram Induction}).
 \xDH
We can easily transform this to logic as follows: We consider ONLY the
nodes
on the paths, e.g. in
Diagram \ref{Diagram Induction} (page \pageref{Diagram Induction}),
part (2.1), we consider the contradictory paths
$x_{0} \xcP x_{1} \xcP x_{2},$ $x_{0} \xcP x_{2},$ and transform them into
$x_{0}= \xCN x_{1},$ $x_{1}= \xCN x_{2},$ $x_{0}= \xCN x_{2}$
and neglect the arrow $x_{1} \xcP x_{3},$ to be added later in step
(2.2.2), and
of course present in Yablo's diagram.
Obviously, this is contradictory in the logical sense.
 \xDH
Consider again
Diagram 
\ref{Diagram Induction} (page 
\pageref{Diagram Induction}), part (2.2.2).
The whole diagram is not contradictory, as we opened a new possibility at
$x_{1},$
which becomes impossible only in part (2.3) of this diagram.
So, as a detail, the paths $x_{0} \xcP x_{1} \xcP x_{2},$ $x_{0} \xcP
x_{2}$ are contradictory, but
the possibility $x_{1} \xcP x_{3}$ (recall that for $x_{0}$ true, $x_{1}$
will be an $'' or''!)$
is closed only in (2.3), so the contradiction of the
whole diagram is only complete in (2.3). If you prefer to think
probabilistically,
in part (2.2.2) $50 \xET $ of the cases are contradictory etc.
See also Section 
\ref{Section C1C2} (page 
\pageref{Section C1C2})  for opening and closing new
possibilities.
 \xEj
 \xEj

\ed

\bd

$\hspace{0.01em}$

% {\tiny (+++ Orig. No.:  Definition Yablo-Structure +++)}

% {\tiny LABEL: {Definition Yablo-Structure}} \\[1mm]
\label{Definition Yablo-Structure}

(Again, see also Section \ref{Section Logic} (page \pageref{Section Logic}).)

 \xEh
 \xDH
(Yablo, see  \cite{Yab82}.)

We call the following graph the Yablo structure, $ \xCf YS:$

Written as graph:

$\{x_{i} \xcP x_{j}:i<j\},$ $i,j \xbe \xbo.$

See Definition \ref{Definition Graph} (page \pageref{Definition Graph}), (3.3))

Written as infinite propositional formulas:

$x_{i}$ $=$ $ \xbf_{x_{i}}$ $=$ $ \xcU \{ \xCN x_{j}:i<j\},$ $i,j \xbe
\xbo.$

Note that we have a very strong coherence property for valuations here:

(R) For $j>i$ $ \xbf_{x_{j}}$ $=$ $ \xbf_{x_{i}} \xex \{x_{k}:k>j\},$
where $ \xex $ stands for ``restriction''.

 \xDH Contradictions in $ \xCf YS:$

Suppose $x_{i}$ is a node in the $ \xCf YS.$ Consider just $x_{i},$
$x_{i+1},$ $x_{i+2},$ and
the arrows between them. (This view is important for the basic
contradictions
in Section 
\ref{Section Elementary-Cells} (page 
\pageref{Section Elementary-Cells})  and
our inductive constructions in
Section 
\ref{Section Combining-Cells} (page 
\pageref{Section Combining-Cells}).)
Then $x_{i} \xcj \xCN x_{i+1} \xcu \xCN x_{i+2},$ and $x_{i+1} \xcj \xCN
x_{i+2},$
which is contradictory.

But starting with $ \xCN x_{i}$ will also result in a contradiction:

$ \xCN x_{i}= \xcO \{x_{j}:i<j\},$ suppose $x_{k}$ is true, then $x_{k}$
is contradictory, as we just
saw, so $ \xCN x_{i}$ is contradictory, too.

Note that we abused notation here, where $ \xcp $ may be an arrow in a
graph,
or an implication in a logical formula. Context will tell.

In more detail:

 \xEh
 \xDH The contradictions have the form $x_{i} \xcP x_{j} \xcP x_{k},$
$x_{i} \xcP x_{k},$ this is the
basic contradiction ``cell'', we call such cells Yablo cells, $YC$'s, or
Yablo triangles.

 \xDH
We sometimes call $x_{i}$ the head, $x_{j}$ the knee, $x_{k}$ the foot of
such $YC$'s.
This is relative to the cell discussed. When we combine such cells,
the knee or foot of one cell may become the head of the next cell, etc.

 \xDH
The distinction between knee and head is not only conceptually
important (at the knee, the contradiction is not complete, at the
head it is complete), the full importance is seen in
Remark 
\ref{Remark Induction-New-1} (page 
\pageref{Remark Induction-New-1}), where we see the fundamental
difference
between continuing the construction at the knee, and at the head.

 \xDH
We sometimes write $x_{i} \xcP_{ \xbD_{i}}x_{i+1}$ for $x_{i} \xcP y_{i}
\xcP x_{i+1},$ $x_{i} \xcP x_{i+1},$ context will
tell. This is used for brevity in
Section 
\ref{Section Combining-Saw-Blade} (page 
\pageref{Section Combining-Saw-Blade}).

 \xDH When the contradiction is not formed by single arrows, but more
generally by
paths, we call this a contradiction cell.

 \xEj

 \xDH We say that a graph has the Yablo Property iff it has a node $x$
such that
both the assumptions $x$ and $ \xCN x$ lead to a contradiction.

 \xEj

\ed

We need infinite depth and width in our constructions:

\bfa

$\hspace{0.01em}$

% {\tiny (+++ Orig. No.:  Fact Infinite +++)}

% {\tiny LABEL: {Fact Infinite}} \\[1mm]
\label{Fact Infinite}

(We claim no originality, and this must be part of the folklore of
the subject.)

 \xEh
 \xDH The construction needs infinite depth,

 \xDH the logic as used in Yablo's construction is not compact,

 \xDH the construction needs infinite width, i.e. the logic cannot be
classical.
 \xEj

\efa

\subparagraph{
Proof
}

$\hspace{0.01em}$

% {\tiny (+++ Orig.:  Proof +++)}

 \xEh
 \xDH
Let $x_{i}$ be a minimal element, then we can chose an arbitrary truth
value
$x_{i},$ and propagate this value upwards. If there are no infinite
descending
chains, we can do this for the whole construction.

 \xDH
The logic as used in Yablo's construction is not compact:

Trivial.
(Take $\{ \xcO \{ \xbf_{i}:i \xbe \xbo \}\} \xcv \{ \xCN \xbf_{i}:i \xbe
\xbo \}.$ This is obviously inconsistent,
but no finite subset is.).

 \xDH
It is impossible to construct a Yablo-like structure with classical logic:

Take an acyclic graph, and interpret it as in Yablo's construction.
Wlog., we may assume the graph is connected. Suppose it shows that
$x_{0}$ cannot be given a truth value. Then the set of formulas showing
this does not have a model, so it is inconsistent. If the formulas were
classical, it would have a finite, inconsistent subset, $ \xbF.$ Define
the depth
of a formula as the shortest path from $x_{0}$ to this formula.
There is a (finite) $n$ such that all formulas in $ \xbF $ have depth $
\xck n.$
Give all formulas of depth $n$ (arbitrary) truth values, and work
upwards using truth functions. As the graph is acyclic, this is possible.
Finally, $x_{0}$ has a truth value.

Thus, we need the infinite $ \xcU / \xcO.$
 \xEj

$ \xcz $
\\[3ex]

We now turn to the third truth value, $ \xbx.$

See also Section7.1.2, Fact 7.1.1 in
 \cite{Sch22}
for a less systematic approach.

\br

$\hspace{0.01em}$

% {\tiny (+++ Orig. No.:  Remark Xi +++)}

% {\tiny LABEL: {Remark Xi}} \\[1mm]
\label{Remark Xi}

Recall the definition of $ \xbx $ in
Definition \ref{Definition Graph} (page \pageref{Definition Graph}), (4.4), and
Section \ref{Section Logic} (page \pageref{Section Logic}).

The truth value $ \xbx $ has the following properties:

 \xEh
 \xDH $ \xCN \xbx = \xbx $ (and, of course $ \xCN \xdT = \xdF,$ $ \xCN
\xdF = \xdT)$
 \xDH the following order $<$ is compatibel with $ \xcu $ and $ \xco:$
$ \xdF < \xbx < \xdT,$ i.e. for any truth values $ \xba,$ $ \xbb $ $
\xba \xcu \xbb =min(\xba, \xbb)$ and
$ \xba \xco \xbb =max(\xba, \xbb).$

This completely describes the behaviour of $ \xdT,$ $ \xbx,$ $ \xdF.$
 \xDH This is the only 3-valued logic with the property $ \xCN \xbx = \xbx
$ and
compatibility with $ \xcu $ and $ \xco.$ It is thus identical to the
Strong Kleene logic
SK, as was pointed out by a referee.
 \xEj

\er

Proof:

 \xEh
 \xDH trivial
 \xDH Consider two copies of $ \xCf YS$'s, $ \xCf YS$ and $ \xCf YS',$
beginning with $x_{0}$ and $x'_{0}$
respectively, and the following two structures
(a) $x \xcP y \xcP z,$ $x \xcP z$ with $x= \xCN y \xcu \xCN z,$ $y= \xCN
z$
(b) $x' \xcP y' \xcP z',$ $x' \xcP z' $ with $x' = \xCN y' \xco \xCN z'
,$ $y' = \xCN z' $

Consider a new origin, $u,$ and arrows from $u$ to $x_{0},$ $x'_{0},$ $x,$
$x' $ as follows:
$u=x_{0} \xcu x'_{0},$ and $u=x_{0} \xcu x,$ $u=x_{0} \xcu x',$
$u=x_{0} \xco x'_{0},$ and $u=x_{0} \xco x,$ $u=x_{0} \xco x'.$

 \xEh
 \xDH
Let $ \xfA \xbf \xfA $ be the truth value of $ \xbf.$

Let $ \xfA \xbf \xfA = \xbx $ be an abbreviation for $ \xfA \xbf \xfA =
\xdF $ and $ \xfA \xCN \xbf \xfA = \xdF.$

Let $ \xfA \xbf \xfA = \xbx.$ We then have

 \xEh
 \xDH $ \xbx \xcu \xdT = \xbx:$
$ \xbf \xcu \xdT = \xdF \xcu \xdT,$ $ \xCN (\xbf \xcu \xdT)= \xCN \xbf
\xco \xdF = \xdF \xco \xdF = \xdF.$
 \xDH $ \xbx \xco \xdT = \xdT:$
$ \xbf \xco \xdT = \xdT,$ $ \xCN (\xbf \xco \xdT)= \xCN \xbf \xcu \xdF
= \xdF \xcu \xdF = \xdF.$
 \xDH $ \xbx \xcu \xdF = \xdF,$ $ \xbx \xco \xdF = \xbx:$
analogously.
 \xDH Likewise $ \xbx \xcu \xbx = \xbx \xco \xbx = \xbx.$
 \xEj

 \xDH We thus have:

$ \xdF \xcu \xbx $ $=$ $ \xCN (\xCN (\xdF \xcu \xbx))$ $=$ $ \xCN (
\xCN \xdF \xco \xCN \xbx)$ $=$ $ \xCN (\xdT \xco \xbx)$ $=$ $ \xCN \xdT
$ $=$ $ \xdF $ and

$ \xdF \xco \xbx $ $=$ $ \xCN (\xCN (\xdF \xco \xbx))$ $=$ $ \xCN (
\xCN \xdF \xcu \xCN \xbx)$ $=$ $ \xCN (\xdT \xcu \xbx)$ $=$ $ \xCN \xbx
$ $=$ $ \xbx $

 \xDH
Using the operations $ \xcu,$ $ \xco $ for truth values $ \xba,$ $ \xbb
$

$ \xba \xco \xbb $ $=$ $sup\{ \xba, \xbb \},$ $ \xba \xcu \xbb $ $=$
$inf\{ \xba, \xbb \}$ and above results,
we have
 \xEh
 \xDH $sup\{ \xbx, \xdF \}= \xbx,$ $sup\{ \xbx, \xdT \}= \xdT $ (or
$inf\{ \xbx, \xdF \}= \xdF,$ $inf\{ \xbx, \xdT \}= \xbx))$ and thus
the
order

 \xDH $ \xdF < \xbx < \xdT,$

 \xDH (and finally the inverse operations

$- \xbx = \xbx,$

$- \xdT = \xdF,$

$- \xdF = \xdT.)$

 \xEj

 \xEj

 \xDH
$ \xdT \xcu \xbx =inf\{ \xdT, \xbx \},$ suppose $ \xdT < \xbx,$ so $
\xCN \xbx = \xbx = \xCN \xdT = \xdF,$ thus $ \xdT < \xbx \xdF,$
contradiction. So $ \xbx < \xdT,$ and $ \xdF = \xCN \xdT < \xCN \xbx =
\xbx.$

 \xEj

\bd

$\hspace{0.01em}$

% {\tiny (+++ Orig. No.:  Definition Path-Value +++)}

% {\tiny LABEL: {Definition Path-Value}} \\[1mm]
\label{Definition Path-Value}

(This applies only to unique occurrences of a variable in the formula
attached to another variable.)

In most cases, $ \xbf_{x}$ will be of the form $ \xcU \xCN x_{i}.$ We
slightly weaken this to
$ \xcU \xCN x_{i}:i \xbe I \xcu \xcU x_{j}:j \xbe J$ - with no
$x_{i}=x_{j}.$ We then can describe $ \xbf_{x}$ in the graph by
$\{x \xcp x_{i}:i \xbe I\} \xcv \{x \xcP x_{j}:j \xbe J\}.$

We can then describe a path $ \xbs $ from $x$ to $y$ by $x \xcp x' \xcP
x''  \Xl y,$ etc., in particular,
it is interesting whether the sign changes by traversing the path:
if the number of $ \xcP $'s is even, it does not, if it is odd, it does.

Thus, we may define a polarity of a path, $pol(\xbs)$ with $pol(\xbs
):=0$ iff it
does not change, and 1 otherwise.

Contradictions between two paths $ \xbs $ and $ \xbs' $ exist iff they
have different
polarity.

Formally:

Let $ \xbs,$ $ \xbs' $ be paths as above.
 \xEh
 \xDH
If $ \xbs:=a \xcp b,$ then $pol(\xbs)=0,$ if $ \xbs:=a \xcP b,$ then
$pol(\xbs)=1.$
 \xDH
Let $ \xbs \xDO \xbs' $ be the concatenation of $ \xbs $ and $ \xbs'.$
Then $pol(\xbs \xDO \xbs')=0$ iff $pol(\xbs)=pol(\xbs'),$ and 1
otherwise.
 \xEj
\subsubsection{
A Remark on More General Formulas
}

% {\tiny LABEL: {Section DNF}} \\[1mm]
\label{Section DNF}

\ed

\br

$\hspace{0.01em}$

% {\tiny (+++ Orig. No.:  Remark DNF-2 +++)}

% {\tiny LABEL: {Remark DNF-2}} \\[1mm]
\label{Remark DNF-2}

Considering more general propositional formulas, beyond $ \xcU x_{i}$ (or
$ \xcU \xCN x_{i})$
seems to entail many complications.

Not surprisingly, they are related to those of the Diamond,
see Section 
\ref{Section Diamond-Problem} (page 
\pageref{Section Diamond-Problem})
(excessive branchings, synchronisation) and
postponed contradictions, see  \cite{Sch25},
version 3, Diagram 4.3 (the HTML version), and
 \cite{Sch23b}, Section 2.

The following is a very simple example:

\er

\be

$\hspace{0.01em}$

% {\tiny (+++ Orig. No.:  Example Ebenen +++)}

% {\tiny LABEL: {Example Ebenen}} \\[1mm]
\label{Example Ebenen}

 \xEh
 \xDH Let
$x_{i}$ $=$ $x'_{i}$ $=$ $(\xcU \xCN x_{j}:j>i)$ $ \xco $ $(\xcU \xCN
x'_{j}:j>i),$ $i \xcg 0,$ so
$ \xCN x_{i}$ $=$ $ \xCN x'_{i}$ $=$ $ \xcO \{x_{j} \xcu x'_{j' }:$ $j,j'
>i\}.$

 \xDH
The structure has the Yablo property:

We argue as for the Yablo structure.
Chose $i,$ and suppose $x_{i}= \xdT,$ so
$x_{i}$ $=$ $(\xcU \xCN x_{j}:j>i)$ $ \xco $ $(\xcU \xCN x'_{j}:j>i).$
Suppose $ \xcU \xCN x_{j}:j>i$ holds, so $ \xCN x_{i+1}$ $=$
$ \xcO \{x_{j} \xcu x'_{j' }:$ $j,j' >i+1\},$ so for some $j>i+1$ $x_{j},$
holds, contradiction.
The case $ \xcU \xCN x'_{j}:j>i$ is analogous.

Chose $i,$ and suppose $x_{i}= \xdF,$ so for some $j,j' >i$ $x_{j} \xcu
x'_{j' }$ holds,
again a contradiction, by the above.

 \xDH Note that the structure has the property (R) of
Definition 
\ref{Definition Yablo-Structure} (page 
\pageref{Definition Yablo-Structure}).

Obviously, the structure may be generalized to more than two disjuncts.
Instead of a single chain, we have a chain of levels.

 \xEj

\ee

\br

$\hspace{0.01em}$

% {\tiny (+++ Orig. No.:  Remark Trivial +++)}

% {\tiny LABEL: {Remark Trivial}} \\[1mm]
\label{Remark Trivial}

Classical formulas in Example 
\ref{Example Ebenen} (page 
\pageref{Example Ebenen})
tend to trivialise the problem.

\er

Let $x_{i} \xco x'_{i}= \xbf_{i} \xco \xbf'_{i},$ $ \xbf'_{i}$ a
conjunction.

If $ \xbf'_{i}= \xdT,$ then $ \xbf_{i} \xco \xbf'_{i}= \xdT.$ If $
\xdT \xbe \xbf'_{i},$ then $ \xbf'_{i}$ $=$ $ \xbf'_{i}$ without $ \xdT
.$

If $ \xbf'_{i}= \xdF,$ then $ \xbf_{i} \xco \xbf'_{i}= \xbf_{i}.$ If $
\xdF \xbe \xbf'_{i},$ then $ \xbf'_{i}= \xdF.$

In particular, finite constructions on one side, which may be
initialised ad libitum, trivialise the whole construction.
\subsection{
Strategy
}

% {\tiny LABEL: {Section Strategy}} \\[1mm]
\label{Section Strategy}

These comments will be clearer in hindsight, after having read the
details.
Still the best place to put them seems to be in the general remarks, and
to go
into details later on.

We show that Yablo's construction is in a certain sense universal:
 \xEh
 \xDH
The basic contradiction cells in Yablo's construction are necessary,
we need three negative sides, neither two nor more than three will do.
See Section 
\ref{Section Elementary-Cells} (page 
\pageref{Section Elementary-Cells}).
 \xDH
Yablo's composition of contradiction cells is nesessary, too. We may do
more,
but a more baroque construction will contain Yablo's construction.
See Section 
\ref{Section Combining-Cells} (page 
\pageref{Section Combining-Cells})  and
``Saw Blade construction'' e.g. in  \cite{Sch22}, for a more baroque
construction.
 \xDH
Consequently, Yablo's construction is sufficient and necessary for some
$x$ such that $x$ and $ \xCN x$ are both contradictory in some
non-circular
directed graph.

Thus, if we find an (image of) Yablo's construction in a non-circular
graph,
this graph will allow some such $x.$
(See Section \ref{Section Injection} (page \pageref{Section Injection})).

Conversely, if a non-circular graph will allow some such $x,$ then we
must be able to find a substructure of the Yablo structure type.
(See Section \ref{Section Converse} (page \pageref{Section Converse})).

 \xEj
\subsubsection{
From Finite to Infinite
}

% {\tiny LABEL: {Section Finite-Infinite}} \\[1mm]
\label{Section Finite-Infinite}

Consider as always graphs with attached formulas of the type $ \xcU
x_{i},$
let the graph be finite.

Fix some node $ \xCf a,$ say e.g. $ \xfA a \xfA = \xdT,$ and let $b$ be a
direct or indirect
successor
of $ \xCf a,$ and let also $ \xfA b \xfA = \xdT.$ Suppose we add at $b$ a
new arrow, say $b \xcp c.$
Then this is a restriction on $b$ (recall, we have a conjunction at $b)$
and $ \xCf a,$
so the set of models is
restricted, as $c$ (or $ \xCN c)$ has to be satisfied. In some cases, we
may even
have a contradiction, e.g. $b$ or even $ \xCf a$ might be impossible.
If, however, $ \xCf \xfA b \xfA = \xdF $ (recall, we now have a
disjunction at $ \xfA b \xfA = \xdF),$ then
we add new possibilities.

Consider now the same $ \xCf a,$ but set $ \xfA a \xfA = \xdF,$ then $b$
will also be $ \xfA b \xfA = \xdF.$
Thus increasing possibilities for $ \xfA a \xfA = \xdT,$ will decrease
possibilities for
$ \xfA a \xfA = \xdF,$ and vice versa.

This is the trivial fact that for a fixed language and classical (finite!)
propositional logic, the set of models is fixed, and every model is
either a model of $ \xbf $ or of $ \xCN \xbf.$ We underline this, as it
is crucial
for our construction of Yablo's structure from finite parts.

Thus, there is no way to approximate the effect of Yablo's
construction by decreasing the set of models of $ \xCf a$ and of $ \xCN
a.$ Decreasing
one, increases the other. The only way is to eliminate one model, at the
price of creating one or more new ones, eliminating the new ones, and so
on.
This is what we do, and there is no other way.

The strategy is outlined in
Section \ref{Section C1C2} (page \pageref{Section C1C2}), and elaborated in
Section \ref{Section Combining-Cells} (page \pageref{Section Combining-Cells}).
\subsubsection{
The Antagonistic Conditions (C1) and (C2)
}

% {\tiny LABEL: {Section C1C2}} \\[1mm]
\label{Section C1C2}

\bd

$\hspace{0.01em}$

% {\tiny (+++ Orig. No.:  Definition C1C2 +++)}

% {\tiny LABEL: {Definition C1C2}} \\[1mm]
\label{Definition C1C2}

Conditions (C1) and (C2)

Consider $x_{0}$ in the Yablo construction (but this applies to all
$x_{i}).$

We denote by (C1) (for $x_{0})$ the condition that $x_{0}= \xdT $ has to
be contradictory,
and by (C2) the condition that $x_{0}= \xdF $ has to be contradictory,
which means
by $ \xCN x_{0}= \xcO \{x_{i}:i>0\}$ that all $x_{i},i>0$ have to be
contradictory, roughly,
in graph language, that all paths from $x_{0}$ have to lead to a
contradiction.

We want to make $x_{0}$ and $ \xCN x_{0}$ contradictory. We know that this
is
impossible in finite graphs, so we have to show how we solve the problem
in the (infinite) limit. We proceed inductively as Yablo did, but we
also see that we basically have no other choice (modulo minor
modifications).

(We need infinite depth and width in the construction, see
e.g.  \cite{Sch25}, Fact 2.4, or
 \cite{Sch23b}, Fact 1.4).

So we want to express by an infinite graph that both
$x_{0}$ and $ \xCN x_{0}$ are inconsistent,
see Definition \ref{Definition Graph} (page \pageref{Definition Graph}).
Therefore, we treat such contradictions as basic building blocks, as Yablo
did.

The problems are then
 \xEh
 \xDH Which kinds of such building blocks exist? See
Section 
\ref{Section Elementary-Cells} (page 
\pageref{Section Elementary-Cells}).
 \xDH How do we combine them to obtain a suitable graph expressing that
both $x$ and $ \xCN x$ are impossible?

The problem is that, if $x$ is contradictory, then $ \xCN x$ is a
tautology
in classical logic, and thus in the theory expressed by a finite acyclic
graph. So we alternate between $x$ and $ \xCN x$ a tautology, and it is
not
clear what such constructions are in the limit, when we leave classical
logic.

This view is too abstract, we have to look at details, as they are given
in the proof of validity for Yablo's construction, or our
reconstruction in
Section \ref{Section Combining-Cells} (page \pageref{Section Combining-Cells}),
we have to look at instances.
 \xEh
 \xDH We start with satisfying (C1) for $x_{0},$ see (2.1) in
Diagram 
\ref{Diagram Induction} (page 
\pageref{Diagram Induction}), $ \xbD_{0}.$ This is necessary, as we
do not postpone
solutions without necessity, and because the formulas attached to
the points have the form $ \xcU.$
 \xDH This creates new problems for (C2), at $x_{1}$ and $x_{2}.$
Consider $x_{1}.$ We satisfy (C2) at $x_{1},$ see (2.2.1) and (2.2.2) in
Diagram 
\ref{Diagram Induction} (page 
\pageref{Diagram Induction}), adding $ \xbD_{1}.$ Again, this is
necessary.
 \xDH But $x_{1}$ is now a branching point. The branch $x_{1}-x_{2}$
continues to be
satisfied for (C1), but the (new) branch $x_{1}-x_{3}$ opens a new problem
for
(C1), as $x_{0}= \xdT $ implies $x_{1}= \xdF,$ so we have an ``OR'' at
$x_{1}= \xdF.$
 \xDH So we attack this new problem in (2.3) in
Diagram \ref{Diagram Induction} (page \pageref{Diagram Induction}),
by adding a new contradiction to this new problem for (C1), etc.
 \xDH So we have new instances of problems for (C1), (C2), we solve the
new instances, say for (C1), the old instances for both (C1) and (C2) stay
solved.

Thus, in the limit, new problems (instances of problems, more precisely)
are
created, they are solved in one of the next steps, problems once solved
stay
solved, so there remain no unsolved problems. All these steps are
necessary.

See the details in
Section 
\ref{Section Combining-Yablo} (page 
\pageref{Section Combining-Yablo})  and
Section 
\ref{Section Combining-Saw-Blade} (page 
\pageref{Section Combining-Saw-Blade}).

The set of problems increases, but all are solved in the limit, and
$x$ and $ \xCN x$ are contradictory.
 \xEj
 \xDH We have (essentially) an alternation of creating and treating
problems
for (C1) and (C2): We cannot treat e.g. a (C1) problem, bevor it was
created by solving a (C2) problem - and conversely.

(Depending on organization, we have some liberty: e.g. solving (C1)
creates a (C2) problem, but this needs TWO branches. So we might first
treat these two branches ((C1) problems) before going back to
the new (C2) problem created by the first new branch.)

 \xEj
\subsubsection{
The Escape Problem
}

\ed

\be

$\hspace{0.01em}$

% {\tiny (+++ Orig. No.:  Example Procrastination +++)}

% {\tiny LABEL: {Example Procrastination}} \\[1mm]
\label{Example Procrastination}

This Example shows that
infinitely many finitely branching points cannot always
replace infinite branching - there is an infinite
``escape branch'' or path.
See below, Definition 
\ref{Definition Escape} (page 
\pageref{Definition Escape}).
This modification of the Yablo structure has one acceptable
valuation for $y_{1}:$

Let $y_{i},$ $i< \xbo $ as usual, and introduce new $x_{i},$ $3 \xck i<
\xbo.$

Let $y_{i} \xcP y_{i+1},$ $y_{i} \xcp x_{i+2},$ $x_{i} \xcP y_{i},$ $x_{i}
\xcp x_{i+1},$ with

$y_{i}:=$ $ \xCN y_{i+1} \xcu x_{i+2},$ $x_{i}:= \xCN y_{i} \xcu x_{i+1}.$

If $y_{1}= \xdT,$ then $ \xCN y_{2} \xcu x_{3},$ by $x_{3},$ $ \xCN y_{3}
\xcu x_{4},$ so, generally,

if $y_{i}= \xdT,$ then $\{ \xCN y_{j}:$ $i<j\}$ and $\{x_{j}:$ $i+1<j\}.$

If $ \xCN y_{1},$ then $y_{2} \xco \xCN x_{3},$ so if $ \xCN x_{3},$
$y_{3} \xco \xCN x_{4},$ etc., so, generally,

if $ \xCN y_{i},$ then $ \xcE j(i<j,$ $y_{j})$ or $ \xcA j\{ \xCN x_{j}:$
$i+1<j\}.$

Suppose now $y_{1}= \xdT,$ then $x_{j}$ for all $2<j,$ and $ \xCN y_{j}$
for all $1<j.$
By $ \xCN y_{2}$ there is $j,$ $2<j,$ and $y_{j},$ a contradiction, or $
\xCN x_{j}$ for all $3<j,$
again a contradiction.

But $ \xCN y_{1}$ is possible, by setting $ \xCN y_{i}$ and $ \xCN x_{i}$
for all $i.$

Thus, replacing infinite branching by an infinite number of finite
branching does not work for the Yablo construction, as we can always chose
the
``escape'' branch, see below.

\ee

\bd

$\hspace{0.01em}$

% {\tiny (+++ Orig. No.:  Definition Escape +++)}

% {\tiny LABEL: {Definition Escape}} \\[1mm]
\label{Definition Escape}

When we try to make a structure contradictory for some value of the
origin, but find a way (a valuation) to avoid or escape a contradiction,
we call such a valuation an escape path, branch, or possibilty. This is
then a
model consistent with the theory of the graph.

This seems an intuitive description how the structure fails to be
contradictory for this valuation.

\ed

$ \xCO $
% {\tiny KS= wab-Anf   FOUND IN END c:-ks-w-t-framart-phil-yablo-wablo-wba-1.mp}
% \\[1mm]

$ \xCO $
% {\tiny KS= wab-Main   FOUND IN START
% c:-ks-w-t-framart-phil-yablo-wablo-wbe-1.mp} \\[1mm]
\clearpage
\section{
Elementary Cells
}

% {\tiny LABEL: {Section Elementary-Cells}} \\[1mm]
\label{Section Elementary-Cells}
\subsection{
Introduction
}

% {\tiny LABEL: {Section Cell-Intro}} \\[1mm]
\label{Section Cell-Intro}

Recall Definition 
\ref{Definition Yablo-Structure} (page 
\pageref{Definition Yablo-Structure})
and Definition \ref{Definition C1C2} (page \pageref{Definition C1C2}).

We try to generalize Yablo cells (i.e. of the type
$x \xcP y \xcP z,$ $y \xcP z,$ with the meaning $x \xcr \xCN y \xcu \xCN
z,$ $y \xcr \xCN z)$
to some other types of cells.

 \xEh
 \xDH
(C1) They will have the form of two paths $ \xbs,$ $ \xbs' $ which
branch say at $x,$ and
meet again at $y$ (but not before), and are contradictory.

(C2) Both branches will have to lead to new contradiction cells if $ \xfA
x \xfA = \xdF $
(or $ \xfA x \xfA = \xbx).$
 \xDH
We require that it is possible to construct a structure $S$ with such
cells
starting,
say at $x_{0}$ such that $ \xfA x_{0} \xfA = \xdT $ and $ \xfA x_{0} \xfA
= \xdF $ are contradictory (classically),
i.e. $ \xfA x_{0} \xfA = \xbx $ using our truth value $ \xbx.$
In particular, we have to pay attention that there are no escape
possibilities,
see Definition \ref{Definition Escape} (page \pageref{Definition Escape}).
 \xDH
In addition, we require that this is possible using only the contradiction
cells of this type. This prevents ``cheating'', where the ``real'' work
is done by some other type of cell. See
Definition 
\ref{Definition Cell-Hierarchy} (page 
\pageref{Definition Cell-Hierarchy}), (1),
and Remark \ref{Remark ESC-Simple} (page \pageref{Remark ESC-Simple}).
 \xDH
We impose syntactical limitations on the formulas attached to
the nodes in the cells see
Section 
\ref{Section Outline} (page 
\pageref{Section Outline}), prerequisites, (2).
 \xDH
These cells will be part of larger structures, so we have to consider
their ``environment'', in particular they must not permit
escape paths by a suitable choice of cases,
see Definition \ref{Definition Escape} (page \pageref{Definition Escape}).

 \xDH
We will see that under these conditions all suitable forms of
contradiction cells are triangles with three negative sides.
This also holds when we consider cells formed by paths, and
not only those formed by single arrows.
 \xEj
\subsection{
Basic Discussion
}

% {\tiny LABEL: {Section Cells-General}} \\[1mm]
\label{Section Cells-General}

\bd

$\hspace{0.01em}$

% {\tiny (+++ Orig. No.:  Definition Cell-Hierarchy +++)}

% {\tiny LABEL: {Definition Cell-Hierarchy}} \\[1mm]
\label{Definition Cell-Hierarchy}

 \xEh
 \xDH The trivial contradiction, the simplest contradiction cell:

$x \xeH y,$ with the meaning $x=y \xcu \xCN y.$

 \xDH The Yablo Cell:

$x \xcP y \xcP z,$ $x \xcP z,$ with the meaning $x= \xCN y \xcu \xCN z,$
$y= \xCN z.$

 \xDH The diamond:

$x \xcP y \xcP z,$ $x \xcP y' \xcp z,$ with the meaning $x= \xCN y \xcu
\xCN y',$ $y= \xCN z,$ $y' =z.$

The diamond will be discussed in
Section \ref{Section Diamond-Problem} (page \pageref{Section Diamond-Problem}).

 \xEj

\ed

We will argue that it suffices to consider these types of cells,
see Remark \ref{Remark Yablo-Xi} (page \pageref{Remark Yablo-Xi})  and
Section \ref{Section PathCells} (page \pageref{Section PathCells}).
Furthermore, we can also neglect all but the (slightly generalized)
Yablo Cells, as we will see in the present section.

It will become clear in a moment that above cells are fundamentally
different, but for this
we have to consider the second requirement (C2).

\br

$\hspace{0.01em}$

% {\tiny (+++ Orig. No.:  Remark ESC-Simple +++)}

% {\tiny LABEL: {Remark ESC-Simple}} \\[1mm]
\label{Remark ESC-Simple}

 \xEh
 \xDH There is an escape problem
(see Definition \ref{Definition Escape} (page \pageref{Definition Escape}))
for $ \xCN x$ in the trivial contradiction,
$x=y \xcu \xCN y,$ so $ \xCN x=y \xco \xCN y.$
(See Definition 
\ref{Definition Cell-Hierarchy} (page 
\pageref{Definition Cell-Hierarchy}).)
If we continue $y$ with $y=z \xcu \xCN z,$ so $ \xCN y=z \xco \xCN z,$
etc. we have an escape
sequence $ \xCN x,$ $ \xCN y,$ $ \xCN z,$ etc., which never meets a
contradiction,
thus providing a model, and failing condition (C2).

 \xDH
Suppose we have an escape problem, like $x=y \xcu \xCN y,$ thus $ \xCN x=y
\xco \xCN y.$

Can we append a new finite structure $ \xbG $ to $y$ which reduces the
possibilities to just one value, say for $z?$ So, whatever the input,
$y$ or $ \xCN y,$ the outcome is $z?$ (``Input'' etc. is seen from the
graph
perspective.)

So we would have $x=y \xcu \xCN y,$ $y= \xbf (z)$ (i.e. a formula which
depends on $z),$ e.g.
$y=z \xcu \xCN z.$

This, however, is not possible.

There are, whatever the inner structure, just four possibilities:
$ \xbf (z)$ might be equivalent to $y=z,$ or to $y= \xCN z,$ or to $y=z
\xcu \xCN z,$ or to $y=z \xco \xCN z.$

The first two cases are trivial, as we have the possibilities $y$ or $
\xCN y,$
so $ \xCN z$ is always a possibility.

Suppose $y=z \xcu \xCN z.$ Then $ \xCN x=y \xco \xCN y=(z \xcu \xCN z)
\xco (\xCN z \xco z)= \xCN z \xco z.$

Suppose $y=z \xco \xCN z.$ Then $ \xCN x=(z \xco \xCN z) \xco \xCN y=z
\xco \xCN z.$

 \xEj

\er

\br

$\hspace{0.01em}$

% {\tiny (+++ Orig. No.:  Remark Simple-Cells +++)}

% {\tiny LABEL: {Remark Simple-Cells}} \\[1mm]
\label{Remark Simple-Cells}

We consider here cells with three arrows.

Note that the following considerations apply also to cells formed by
paths consisting of more than one arrow, and not only to cells formed
by single arrows.
See Section \ref{Section PathCells} (page \pageref{Section PathCells}).
In particular, also in more complicated cells all sides have to negative.

In the graph with not annotated arrows, they have the form
$x-y-z$, $x-z$, where the lines may stand for $ \xcp $ or $ \xcP.$

Again, we want a contradiction for $x$ positive, so we need an $ \xcu $ at
$x.$
For the same reason, the number of negative arrows may not be even.

By (C2), $y$ and $z$ must be negative, when $x$ is positive, so
$x-y$ and $x-z$ must be negative, $x \xcP y$ and $x \xcP z.$ Thus, all
three must be
negative, so $y \xcP z$ too.

This is the original type of contradiction in Yablo's construction

$x \xcP y \xcP z,$ $x \xcP z.$

This will be discussed in detail in the rest of the paper.
But we see already that both paths, $x \xcP y$ and $x \xcP z$ change sign,
so $y$ and $z$ will be positive when $x$ is negative, appending the same
type of cell
at $y$ and $z$ solves the problem (locally), and offers no escape.

\er

See also the discussion of ``Saw Blades''
in  \cite{Sch22}, Section 3.

 \cite{Sch25}, in particular $[v7],$ the proof of
Remark 2.2 there presents more arguments showing that all arrows have
to be negative.
\subsection{
Values for y and z in Yablo Cells
}

% {\tiny LABEL: {Section YC-Details}} \\[1mm]
\label{Section YC-Details}

We consider here what happens when we give $z$ and $y$ truth values $ \xba
$ $(\xdT,$ $ \xdF,$ $ \xbx).$
In addition, we consider the variants $y= \xCN z \xcu \xba $ and $y= \xCN
z \xco \xba.$

\br

$\hspace{0.01em}$

% {\tiny (+++ Orig. No.:  Remark Yablo-Xi +++)}

% {\tiny LABEL: {Remark Yablo-Xi}} \\[1mm]
\label{Remark Yablo-Xi}

See Remark \ref{Remark Xi} (page \pageref{Remark Xi}).

We present a systematic treatment of variants of the ``environment''
of the Yablo triangle.

 \xEh
 \xDH In the Yablo construction, we attach to the Yablo
triangle $x \xcP y \xcP z,$ $x \xcP z,$ to $y$ and $z$
constructions, which make $y= \xCN z \xcu \xbx,$ $z= \xbx.$
We also consider here the cases where we end $z$ by $ \xdT $ or $ \xdF,$
and set e.g. $y= \xCN z \xcu \xdT,$ $y= \xCN z \xco \xdT,$ etc., see
(3.1) below.

Note that $x= \xCN z \xcu \xCN y$ will always hold.

 \xDH We consider the following requirements:
 \xEh
 \xDH Do we have a contradiction for $x= \xdT,$ i.e. will $x$ not be true
in the
Yablo Cell?
 \xDH Will both branches lead to possible contradictions for $x= \xdF?$
 \xDH Do we obtain $x= \xbx?$
 \xEj

 \xDH We consider the cases in all possible combinations:

 \xEh
 \xDH

(a)

$ \xBc 1.b \xBe $ $z= \xdT,$ $ \xBc 2.b \xBe $ $z= \xdF,$ $ \xBc 3.b \xBe $ $z=
\xbx $

(b)

$ \xBc a.1 \xBe $ $y= \xCN z \xcu \xdT = \xCN z,$ $ \xBc a.2 \xBe $ $y= \xCN z
\xcu \xdF =
\xdF,$ $ \xBc a.3 \xBe $ $y= \xCN z \xcu \xbx,$

$ \xBc a.4 \xBe $ $y= \xCN z \xco \xdT = \xdT,$ $ \xBc a.5 \xBe $ $y= \xCN z
\xco \xdF = \xCN
z,$ $ \xBc a.6 \xBe $ $y= \xCN z \xco \xbx.$

Case $ \xBc a.1 \xBe $ is equivalent to case $ \xBc a.5 \xBe.$

Thus, we look at e.g. case $ \xBc 1.1 \xBe,$ i.e. $z= \xdT $ and $y= \xCN z
\xcu
\xdT = \xCN z,$ etc., and
consider whether the requirements in (2) above are satisfied.

 \xDH The following cases satisfy all requirements in (2):
$ \xBc 2.3 \xBe,$ $ \xBc 3.1 \xBe,$ $ \xBc 3.2 \xBe,$ $ \xBc 3.3 \xBe,$ $
\xBc 3.5 \xBe,$ $ \xBc 3.6 \xBe.$
The proof is a tedious exercise and left to the reader.

 \xDH We now look at the results of embedding these triangles in a larger
structure.

Here, $y_{ \xbx }$ etc. stand for $y$ which is of type $ \xbx,$ or for
some $y' $ appended
at $y$ (context will tell),
as we will append in the full structure at $y$ and $z$ a construction
with value $ \xbx.$

In all cases $x= \xCN y \xcu \xCN z$

\xEn
 \xDH $ \xBc 2.3 \xBe $
$z= \xdF,$ $y= \xCN z \xcu y_{ \xbx }.$ so $x= \xCN y_{ \xbx }.$

 \xDH $ \xBc 3.1 \xBe $
$z=z_{ \xbx },$ $y= \xCN z \xcu \xdT,$ so $x=z_{ \xbx } \xcu \xCN z_{
\xbx }.$

 \xDH $ \xBc 3.2 \xBe $
$z=z_{ \xbx },$ $y= \xdF,$ so $x=z_{ \xbx }.$
 \xDH $ \xBc 3.3 \xBe $ (Yablo)
$z=z_{ \xbx },$ $y= \xCN z \xcu y_{ \xbx }= \xCN z_{ \xbx } \xcu y_{ \xbx
},$ so $ \xCN y=z_{ \xbx } \xco \xCN y_{ \xbx },$ and
$x= \xCN y_{ \xbx } \xcu \xCN z_{ \xbx }.$

 \xDH $ \xBc 3.5 \xBe,$ identical to $ \xBc 3.1 \xBe.$

 \xDH $ \xBc 3.6 \xBe $
$z=z_{ \xbx },$ $y= \xCN z \xco y_{ \xbx }.$ Thus, $ \xCN y=z_{ \xbx }
\xcu \xCN y_{ \xbx },$ $x=z_{ \xbx } \xcu \xCN y_{ \xbx } \xcu \xCN z_{
\xbx }.$
\xEp

 \xEj

 \xDH
Summary of the (in this context) important constructions
 \xEh
 \xDH
The Yablo Cell is of the type $x= \xCN y \xcu \xCN z,$ $y= \xCN z,$ with
$y=z= \xbx,$ type $ \xBc 3.3 \xBe $
above, (see Definition 
\ref{Definition Cell-Hierarchy} (page 
\pageref{Definition Cell-Hierarchy})).

 \xDH
The simplified Saw Blade tooth is of the type $x= \xCN y \xcu \xCN z,$ $y=
\xCN z,$ with $z= \xbx,$
$y= \xdF,$ type $ \xBc 3.2 \xBe $ above, see Section 7.6 in  \cite{Sch22},
or Section 3 in  \cite{Sch23b}.
(In this construction, $y$ will be (except for the start) again of type $
\xbx,$
as we obtain this by transitivity, forced by condition (C1), see
Construction 
\ref{Construction Details} (page 
\pageref{Construction Details}), and
Diagram \ref{Diagram Induction} (page \pageref{Diagram Induction}),
steps (2.3), (2.5), (2.6)).
Both the original Saw Blade and the simplified one contain the Yablo
construction.

See, e.g., Diagram 3.1, upper part, in  \cite{Sch23b}:
$x_{ \xbs,0}$ has the role of $x_{0}$ in Yablo's construction. $x_{ \xbs
,0} \xcP x_{ \xbs,1} \xcP y_{ \xbs,0},$
$x_{ \xbs,0} \xcP y_{ \xbs,0}$ is a Yablo cell. $x_{ \xbs,1}$ must be
the head of a new Yablo cell,
starting e.g. with
$x_{ \xbs,1} \xcP x_{ \xbs,2}.$ $x_{ \xbs,0} \xcP x_{ \xbs,1} \xcP x_{
\xbs,2}$ has to be contradicted, so we need
$x_{ \xbs,0} \xcP x_{ \xbs,2},$ so $x_{ \xbs,2}$ must be the head of a
new Yablo cell, so we have some
$x_{ \xbs,2} \xcP x_{ \xbs,3},$ etc.

 \xEj

 \xDH We obtain infinite width and depth by suitable combinations
of cells of the same type,
see Section \ref{Section Cell-Intro} (page \pageref{Section Cell-Intro}), and
 \cite{Sch25}, Fact 2.4 or
 \cite{Sch23b}, Fact 1.4.
See also Section 
\ref{Section Combining-Cells} (page 
\pageref{Section Combining-Cells}).
This requires continuation at both $y$ and $z$ with constructions of value
$ \xbx,$
case $ \xBc 3.3 \xBe.$

Recall that we also use cells of the type $ \xBc 3.2 \xBe $ for the simplified
Saw
Blade construction. However, the conditions (C1) and (C2) force
transitivity and resulting cells of type $ \xBc 3.3 \xBe,$ see
Section \ref{Section Combining-Cells} (page \pageref{Section Combining-Cells}).

 \xEj
\subsection{
Yablo Cells Formed By Paths
}

% {\tiny LABEL: {Section PathCells}} \\[1mm]
\label{Section PathCells}

\er

We consider a triangle $x- \Xl -y- \Xl -z$, $x- \Xl -z$, where $x- \Xl
-y$ etc. may be
composite paths.

\br

$\hspace{0.01em}$

% {\tiny (+++ Orig. No.:  Remark NegativePaths +++)}

% {\tiny LABEL: {Remark NegativePaths}} \\[1mm]
\label{Remark NegativePaths}

We know from
Remark 
\ref{Remark Simple-Cells} (page 
\pageref{Remark Simple-Cells})  that the paths have to be negative,
but they may
be composed of positive and negative arrows.
More precisely,
Remark \ref{Remark Simple-Cells} (page \pageref{Remark Simple-Cells}), (2)
is trivially true for paths, too, and so is (3.1).

We require that $x= \xdT $ is of type $ \xcU,$ and that the triangle is
contradictory for the
case $x= \xdT,$ as well as for $x= \xdF,$ both paths
$x- \Xl -y- \Xl -z$ und $x- \Xl -z$ have to lead to a contradiction
(conditions (C1) and
(C2)). Note that $y$ and $z$ have to be the first nodes on above paths to
be contradictory for $x= \xdF,$ - otherwise we will not reach $y$ (or
$z),$ the path is
barred before.
\subsubsection{
Branching points
}

% {\tiny LABEL: {Section BranchingPoints}} \\[1mm]
\label{Section BranchingPoints}

\er

\br

$\hspace{0.01em}$

% {\tiny (+++ Orig. No.:  Remark Branching +++)}

% {\tiny LABEL: {Remark Branching}} \\[1mm]
\label{Remark Branching}

 \xEh
 \xDH $y$ is the origin of a new cell, so we have to branch at $y.$
(Anyway, we have to branch on the path $x-y-z$, otherwise, we have the
trivial
contradiction (and an escape possibility).)

Branching on $x-z$ at some additional intermediate point $x' $ does not
necessarily
create a Diamond,
see Definition 
\ref{Definition Cell-Hierarchy} (page 
\pageref{Definition Cell-Hierarchy}),
as we require a new cell only at $y$ and $z,$ and not at $x'.$
Still, this is different from additional branching on $x-y-z$, see
(3) below.

Independence of the choices makes this case uninteresting:

We have to avoid that one choice at $y,$ and one choice at $y' $ do not
lead to a
contradiction. But then all combinations of choices at $y$ and choices at
$y' $
have to lead to contradictions. (Simple argument with quantifiers.)

The argument in Section 
\ref{Section Combining-Cells} (page 
\pageref{Section Combining-Cells})
that all suitable constructions contain Yablo's constructions is valid,
as we have only one knee, and the part head-knee-foot is not affected.

 \xDH Branching at $y:$

Suppose we branch at $y,$ so we have, in addition to the triangle, some
$y-z' $

We know that $x- \Xl -y$ is negative, so, if $ \xfA x \xfA = \xdT,$ then
$ \xfA y \xfA = \xdF.$
This is the situation of the Yablo construction,
see Section 
\ref{Section Combining-Cells} (page 
\pageref{Section Combining-Cells}).
So just a short comment here.

Thus, we also need for $x-y-z' $ a contradiction, this leads to infinite
branching at $x.$

 \xEh
 \xDH A new path $x-z' $ as in Yablo's construction leads to infinite
branching and non-classical logic.

 \xDH We may branch on the already existing $x-y$ oder $ \xCf x-z$, and
continue to
build a (finite) contradiction to $x-y-z'.$

This however, is possible only a finite number of times, e.g.
$x-a_{0}-a_{1}-z',$ then $x-a_{0}-a_{1}-a_{2}-z'',$ etc.,
see Example 
\ref{Example Procrastination} (page 
\pageref{Example Procrastination}).
This constructs an infinite sequence of choices
$a_{0}-a_{1}-a_{2}- \Xl $ which offers an escape possibility.

Thus we need infinite branching (and non-classical logic).

 \xDH
We see that we may consider only the new paths $x-z',$ they form an
infinite subsequence of the construction. By full transitivity of the
final
construction, any infinite subsequence of the construction has the same
properties as the full construction.
 \xEj
 \xDH Branching at $x' $ on the path $x-x' -z$:

Consider the case $x= \xdT.$ Whatever the branching situation on $x-y-z$
is
(just one branching, or multiple branchings), for the branch considered
in $x-y-z$, the path $x-x' -z$ is supposed to create a contradiction.
Suppose we have a branching at $x',$ where $ \xfA x' \xfA = \xdF,$ i.e.
it is a logical ``OR'',
then we may decide at $x' $ otherwise, and not have a contradiction
with the branch $x-y-z$. Thus, we first chose a path via $y,$ and then
chose a contradicting path $x-x' -z$ (or directly $x-z$). In Yablo's
construction, this is trivial, as $x-z$ is a direct link.

See Construction 
\ref{Construction Details} (page 
\pageref{Construction Details}), (3), ``Note''
where this affects the recursive construction.

So we stipulate that, for $ \xfA x \xfA = \xdT,$ any branching point $x'
$ on $x-x' -z$ must
not be $ \xfA x' \xfA = \xdF,$ i.e. not an ``OR''.
 \xEj
\subsubsection{
Paths
}

% {\tiny LABEL: {Section Yablo-Paths}} \\[1mm]
\label{Section Yablo-Paths}

\er

 \xEh
 \xDH Additional branching on $x- \Xl -y$, e.g. $x-x'' -y-z$, $x-z$,
$x'' -y' -z' $

 \xEh
 \xDH Case 1: if $ \xfA x \xfA = \xdT,$ then $x'' $ is a logical OR

This generates a copy of the structure for $ \xfA x \xfA = \xdT:$

We construct for $x-x'' -y' $ as for $x-x'' -y$, i.e. $y' -z',$ $x-z'.$

Note: $y' $ need not be contradictory for $ \xfA x \xfA = \xdF $, as $y$
is contradictory, thus,
$x'' $ will be contradictory ($ \xCf x'' $ will be a logical AND for $
\xfA x \xfA = \xdF),$
so (C2) is satisfied.
Thus, $x'' $ has the role of $y,$ and we may invoke our rule that
contradictions should happen immediately, and not be postponed.

 \xDH Case 2: if $ \xfA x \xfA = \xdT,$ then $x'' $ is a logical AND

So $ \xfA x \xfA = \xdF $ implies $x'' $ is a logical OR,
so we make a copy of the structure for $ \xfA x \xfA = \xdF,$ consider
e.g.
$x-z$, $x-z',$ $x-x'' -y-z$, $x'' -y' -z' $

Note: We simplify, suppose e.g. that we have $x \xcp x'',$ then we do not
need
another contradiction at $x'',$ the one at $x$ suffices (condition (C1)).
We need, however, that all branches lead to contradictions (condition
(C2)).

We might identify $x$ and $x'' $, and impose this as an additional
requirement
for the structures we examine, along
with ``no postponement of contradictions''.

 \xEj

 \xDH We branch on $y- \Xl -z$, and have $x-y-y'' -z$, $y'' -z' $.
 \xEh
 \xDH Case 1: if $ \xfA x \xfA = \xdT,$ then $y'' $ is a logical OR:

We make a copy of the structure for $ \xfA x \xfA = \xdT.$

 \xDH Case 2: if $ \xfA x \xfA = \xdT,$ then $y'' $ is a logical AND:

So $ \xfA x \xfA = \xdF $ implies $y'' $ is a logical OR:

We make a copy of the structure for $ \xfA x \xfA = \xdF.$
 \xEj

Consider e.g. $x-z$, $x-y-y'' -z$, $y'' -z',$ $x-z' $.

 \xDH As branching at $y$ is not different from branching before or after
$y,$ we assume for simplicity that we branch at $y.$ Moreover, as
argued in  \cite{Sch23b},
Section 3.4.1 (``pipeline''), we may assume for simplicity that we always
start contradictions at $x.$
 \xDH If possible, we may also trivialize unwanted branchings as done in
Section \ref{Section Result} (page \pageref{Section Result}).
 \xEj

\br

$\hspace{0.01em}$

% {\tiny (+++ Orig. No.:  Remark Paths +++)}

% {\tiny LABEL: {Remark Paths}} \\[1mm]
\label{Remark Paths}

There are two ways to treat this. We may construct the additional
branchings
the same way as the basic Yablo construction, or we may embed the basic
Yablo construction in the additional branchings as we embed the Yablo
construction in other graphs, see
Section \ref{Section Result} (page \pageref{Section Result}).
\subsection{
The Problem with Diamonds
}

% {\tiny LABEL: {Section Diamond-Problem}} \\[1mm]
\label{Section Diamond-Problem}

\er

\br

$\hspace{0.01em}$

% {\tiny (+++ Orig. No.:  Remark Rhombus-Basic +++)}

% {\tiny LABEL: {Remark Rhombus-Basic}} \\[1mm]
\label{Remark Rhombus-Basic}

See Diagram \ref{Diagram Diamond} (page \pageref{Diagram Diamond}).

 \xEh
 \xDH We have a conflict between the diamonds starting at $y$ and $y' $
and
the diamond starting at $x.$

If $ \xfA x \xfA = \xdT,$ then $ \xfA y \xfA = \xdF $ and $ \xfA y' \xfA
= \xdF $. As the choices at $y$ and $y' $ are
independent, any branch $x-y-y_{1}-z$, $x-y-y_{2}-z$, $x-y-z$
combined with any branch $x-y' -y'_{1}-z$, $x-y' -y'_{2}-z$, $x-y' -z$
must be conflicting, thus, given $x-y-z$ is negative, all branches
on the left must be negative, likewise, all branches on the right must
positive.

However, if $y$ is positive, the diamond $y-y_{1}-z$, $y-y_{2}-z$
has to be contradictory, so not both branches may be negative,
likewise for $y'.$

 \xDH A solution is to ``synchronise'' the choices at $y$ and $y'.$
Informally, this is possible by:

``$ \xCN y \xcu \xCN y',$ and at $y$ choose the arrows $y \xcP y_{2}$ or
$y \xcp z,$ and at $y',$ choose
the arrows $y' \xcP z$ or $y' \xcP y'_{1}$''.

However we express this, the negation is

``$y \xco y' $ or choose at $y$ $y \xcP y_{1}$ or choose at $y' $ the
arrow $y' \xcP y'_{2}$'',

which has nothing to do with the usual negations in Yablo structures.
A way out might be to put the choices at $y$ and $y' $ into a background
theory which is not affected by negation. This is excluded by
Section \ref{Section Outline} (page \pageref{Section Outline}), (1.6).

Even if we consider arbitrary formulas, we cannot speak about choices
at $y$ and $y',$ or arrows like $y \xcP y_{1},$ $y \xcP y_{2}.$
Thus, within our framework, the diamond is NOT a suitable contradiction
cell.

(Yablo's construction does not need synchronisation, this is
done automatically by the universal quantifier.)

 \xDH The simplified version with recursion on one side only, here on the
left
only, is nothing but the Yablo triangle:
$x= \xCN y \xcu \xCN y',$ $y= \xCN z,$ $y' =z \xcu (y'' \xco \xCN y''
\xco \xdT)=z.$

 \xEj

\clearpage

\begin{diagram}

Nested Diamonds, Details

% {\tiny LABEL: {Diagram Diamond}} \\[1mm]
\label{Diagram Diamond}
\index{Diagram Diamond}

\unitlength.8mm
\begin{picture}(280,280)(0,0)

\put(20,270){Example for synchronisation, see Remark \ref{Remark
Rhombus-Basic}}

\put(80,90){\circle*{1}}
\put(78,85){\footnotesize $x$}

\put(80,90){\line(-1,2){40}}
\put(80,90){\line(1,2){40}}

\put(62,130){$-$}
\put(95,130){$-$}

\put(40,170){\circle*{1}}
\put(120,170){\circle*{1}}

\put(30,170){\footnotesize $y$}
\put(122,170){\footnotesize $y'$}

\put(40,170){\line(-1,1){20}}
\put(40,170){\line(0,1){20}}
\put(40,170){\line(1,2){40}}
\put(120,170){\line(1,1){20}}
\put(120,170){\line(0,1){20}}
\put(120,170){\line(-1,2){40}}

\put(25,180){$-$}
\put(36,180){$-$}
\put(122,180){$-$}
\put(133,180){$-$}

\put(20,190){\circle*{1}}
\put(40,190){\circle*{1}}
\put(120,190){\circle*{1}}
\put(140,190){\circle*{1}}

\put(10,190){\footnotesize $y_{1}$}
\put(30,190){\footnotesize $y_{2}$}
\put(121,190){\footnotesize $y'_{1}$}
\put(142,190){\footnotesize $y'_{2}$}

\put(20,190){\line(1,1){60}}
\put(40,190){\line(2,3){40}}
\put(140,190){\line(-1,1){60}}
\put(120,190){\line(-2,3){40}}

\put(46,220){$+$}
\put(57,220){$-$}
\put(69,220){$+$}
\put(89,220){$-$}
\put(101,220){$+$}
\put(112,220){$-$}

\put(80,250){\circle*{1}}
\put(78,255){\footnotesize $z$}

\put(120,170){\line(-2,-1){3}}
\put(120,170){\line(1,-2){1.5}}
\put(40,170){\line(2,-1){3}}
\put(40,170){\line(-1,-2){1.5}}

\put(120,190){\line(-1,-1){2}}
\put(120,190){\line(1,-1){2}}
\put(40,190){\line(-1,-1){2}}
\put(40,190){\line(1,-1){2}}

\put(140,190){\line(-1,0){2}}
\put(140,190){\line(0,-1){2}}
\put(20,190){\line(1,0){2}}
\put(20,190){\line(0,-1){2}}

\put(81,249){\line(0,-1){4}}
\put(81,249){\line(4,-1){4}}
\put(79,249){\line(0,-1){4}}
\put(79,249){\line(-4,-1){4}}

\end{picture}

\end{diagram}

\vspace{10mm}

\clearpage

\subsection{
Conclusion
}

% {\tiny LABEL: {Section Con-2}} \\[1mm]
\label{Section Con-2}

\er

We have shown here that - within our framework - the only useful
contradictions are Yablo cells
(see Definition 
\ref{Definition Yablo-Structure} (page 
\pageref{Definition Yablo-Structure})),
perhaps with paths instead of single arrows.
Simpler contradictions
(see Definition 
\ref{Definition Cell-Hierarchy} (page 
\pageref{Definition Cell-Hierarchy}))
have escape possibilities
(see Definition \ref{Definition Escape} (page \pageref{Definition Escape})),
more complicated ones have synchronisation problems,
see Section 
\ref{Section Diamond-Problem} (page 
\pageref{Section Diamond-Problem}).
\clearpage
\section{
Combining Cells
}

% {\tiny LABEL: {Section Combining-Cells}} \\[1mm]
\label{Section Combining-Cells}
\subsection{
The Main Idea
}

% {\tiny LABEL: {Section Combining-Main}} \\[1mm]
\label{Section Combining-Main}

This section is the core of the present article.

It shows that Yablo's construction is necessary for an acyclic graph
to contain an element $x$ with $x$ and $ \xCN x$ leading to a
contradiction in
the following sense: all other such constructions will contain Yablo's
construction.

We see that there are basically 4 ways to construct such a graph. One
potential possibility (construction only via the head, see (2.2) below
in Remark \ref{Remark Contradiction} (page \pageref{Remark Contradiction})
fails, and all others will contain Yablo's construction as a substructure.

This is necessary, whether we construct such a structure from scratch,
or we find it as a substructure of an existing structure. The result
in Section \ref{Section Result} (page \pageref{Section Result})
is then trivial.
\subsection{
Introduction
}

% {\tiny LABEL: {Section Combining-Intro}} \\[1mm]
\label{Section Combining-Intro}

We refer the reader to
Section 
\ref{Section C1C2} (page 
\pageref{Section C1C2})  for an abstract discussion of Yablo's and
similar
constructions.

We first discuss the different possibilities of composing Contradiction
Cells
to a structure with the Yablo property.
Familiarity with Yablo's construction and the concepts of knee, head, etc.
is
assumed.

See Definition 
\ref{Definition Yablo-Structure} (page 
\pageref{Definition Yablo-Structure}),
and in general
Section \ref{Section Introduction} (page \pageref{Section Introduction})  and
Section 
\ref{Section Elementary-Cells} (page 
\pageref{Section Elementary-Cells}).

It is essential that we need the Contradiction Cells and the way we
compose them in an inductive construction, thus we obtain necessary
constructions. (The sufficiency was shown by Yablo.) We may take
the Cells (and their components) from scratch, or find them in an
existing graph. We need them, and once we hace added a cell or an arrow,
we
need other cells and arrows, by the conditions (C1) and (C2)
(see Definition \ref{Definition C1C2} (page \pageref{Definition C1C2})).
This is the essential property.

\br

$\hspace{0.01em}$

% {\tiny (+++ Orig. No.:  Remark Contradiction +++)}

% {\tiny LABEL: {Remark Contradiction}} \\[1mm]
\label{Remark Contradiction}

Note that we build contradictions. In the case of simple arrows, working
from scratch, this
is easy as we may always chose negative arrows. In the case of paths,
we must be able to choose paths of the right polarity, see
Definition \ref{Definition Path-Value} (page \pageref{Definition Path-Value})
for definition, and
Section \ref{Section Injection} (page \pageref{Section Injection}),
Condition \ref{Condition Injection} (page \pageref{Condition Injection}), and
Example \ref{Example Pol} (page \pageref{Example Pol}).

\er

 \xEh

 \xDH We have 4 possibilities in the inductive construction:
 \xEh
 \xDH we attach a new cell to the
 \xEh
 \xDH head

or the

 \xDH knee
 \xEj

of an old cell,

and

 \xDH if we attach to the knee of an old cell,

 \xEh

 \xDH we may re-use the arrow from knee to head of the old cell to build
the
new cell,

or

 \xDH start from scratch with two new arrows attached to the knee of the
old
cell.

 \xEj

 \xEj

 \xDH For this, we discuss two constructions:
 \xEh
 \xDH The first variant,
see Section 
\ref{Section Combining-Yablo} (page 
\pageref{Section Combining-Yablo}),
is simpler, and Yablo's original construction.

It is the above possibility (1.2.1) ``knee and re-use'',
see Section 
\ref{Section Combining-Yablo} (page 
\pageref{Section Combining-Yablo}).

 \xDH The second variant,
see Section 
\ref{Section Combining-Saw-Blade} (page 
\pageref{Section Combining-Saw-Blade}),
is more complicated (it is the simplified version
of the author's Saw Blade construction as discussed in
 \cite{Sch21} v2 and v3,
 \cite{Sch22},
 \cite{Sch23b}, etc.).

We show in the discussion of the second variant why we have to attach at
the
knee, above (1.1.2), attaching at the head, above (1.1.1), does not
suffice.
(See the Upper Part in Diagram 
\ref{Diagram Induction-New} (page 
\pageref{Diagram Induction-New})).
So this settles the first of above questions.

The second variant is not, but still contains, Yablo's original
construction.
This is a consequence of the interplay of the necessary (and sufficient)
conditions (C1) and (C2),
see Section \ref{Section Strategy} (page \pageref{Section Strategy}).

\br

$\hspace{0.01em}$

% {\tiny (+++ Orig. No.:  Remark Kernel +++)}

% {\tiny LABEL: {Remark Kernel}} \\[1mm]
\label{Remark Kernel}

The $y_{i}$ in
Remark \ref{Remark Induction-New-1} (page \pageref{Remark Induction-New-1})
and Section 
\ref{Section Simplifications} (page 
\pageref{Section Simplifications})
form a kernel
(see e.g.  \cite{Gri26}, $p.$ 331),
thus we may have a kernel in graphs with the Yablo property of
arbitrary cardinality (0 in Yablo's construction, countable in
the simplified Saw Blade construction, of even finite, see e.g.
 \cite{Sch21}, v3, or
 \cite{Sch22},
and if we take an arbitrary set of such construction, also of
uncountable cardinality.

\er

 \xEj

So, whatever we do, we have a structure which contains the original
Yablo construction.

 \xEj

\br

$\hspace{0.01em}$

% {\tiny (+++ Orig. No.:  Remark Disorder +++)}

% {\tiny LABEL: {Remark Disorder}} \\[1mm]
\label{Remark Disorder}

The constructions need not follow an inductive order. We might start with
arbitrary chunks of a graph, the conditions (C1) und (C2) will connect
the arrows via the knees.
\clearpage
\subsection{
Variant Yablo
}

% {\tiny LABEL: {Section Combining-Yablo}} \\[1mm]
\label{Section Combining-Yablo}

\er

\bcs

$\hspace{0.01em}$

% {\tiny (+++ Orig. No.:  Construction Details +++)}

% {\tiny LABEL: {Construction Details}} \\[1mm]
\label{Construction Details}

We now consider the inductive construction of the Yablo structure,
see Diagram \ref{Diagram Induction} (page \pageref{Diagram Induction}).
(The numbers (2.x) in the diagram refer to the points below.)

 \xEh
 \xDH The principle
 \xEh
 \xDH If $x$ has to be the start of a Yablo cell, but is a terminal node,
add
an arrow $x \xcP y.$
 \xDH If $x$ has to be the start of a Yablo cell, $x \xcP y$ exists, but
not $x \xcP z,$
add $x \xcP z$ (transitivity).

If $x=x_{0},$ $z$ will be the start of a new Yablo cell, if $x= \xdF.$
 \xEj
 \xDH Details
 \xEh

 \xDH $x_{0}= \xCN x_{1} \xcu \xCN x_{2},$ $x_{1}= \xCN x_{2},$ thus
$x_{0}= \xdF.$

$x_{1},$ $x_{2}$ will be the start of new Yablo cells.

 \xDH By (C2),
$x_{0} \xcP x_{1}$ has to lead to a contradiction for $x_{0}= \xdF.$

Add $ \xbD_{1}$ by (C2) for $x_{0} \xcP x_{1}.$

 \xEh
 \xDH preparation: we use $x_{1} \xcP x_{2}$ as part of the new Yablo
cell, and add
the new arrow $x_{2} \xcP x_{3}$ as second part of the new cell,
where $x_{3}$ is a new point.

 \xDH finish: We add the new arrow $x_{1} \xcP x_{3}$ as third part of the
new cell.

Note:
$x_{1}= \xCN x_{2}$ was (part of) a full contradiction, this is now a
partial
contradiction, as the new possibility $x_{1}= \xCN x_{3}$ is added.
 \xEj

 \xDH By (C1), repair $ \xbD_{0}:$

$x_{0}= \xCN x_{1} \xcu \xCN x_{2} \xcu \xCN x_{3},$ $x_{1}= \xCN x_{2}
\xcu \xCN x_{3},$ $x_{2}= \xCN x_{3}$

(branch $x_{1} \xcP x_{3}$ has to be contradicted for $x_{0}= \xdT,$ add
$x_{0} \xcP x_{3}),$ thus

$x_{2}= \xCN x_{3},$ $x_{1}= \xdF,$ $x_{0}$ $=$ $ \xdF.$

$x_{3}$ will be the start of a new Yablo cell.

 \xDH By (C2),
$x_{0} \xcP x_{2}$ has to lead to a contradiction for $x_{0}= \xdF $

Add $ \xbD_{2}$ by (C2) for $x_{0} \xcP x_{2}.$

 \xEh
 \xDH preparation: we use $x_{2} \xcP x_{3}$ as part of the new Yablo
cell, and add
the new arrow $x_{3} \xcP x_{4}$ as second part of the new cell,
where $x_{4}$ is a new point.
 \xDH finish: We add the new arrow $x_{2} \xcP x_{4}$ as third part of the
new cell.
 \xEj

 \xDH By (C1), repair $ \xbD_{1}$ (but damage $ \xbD_{0}):$

$x_{0}= \xCN x_{1} \xcu \xCN x_{2} \xcu \xCN x_{3},$ $x_{1}= \xCN x_{2}
\xcu \xCN x_{3} \xcu \xCN x_{4},$ $x_{2}= \xCN x_{3} \xcu \xCN x_{4},$
$x_{3}= \xCN x_{4}$

(branch $x_{2} \xcP x_{4}$ has to be contradicted for $x_{1}= \xdT,$ add
$x_{1} \xcP x_{4}).$ Thus

$x_{3}= \xCN x_{4},$ $x_{2}= \xdF,$ $x_{1}$ $=$ $ \xdF,$ $x_{0}$ $=$
$x_{4}.$

Complete Yablo cell starting at $x_{1}.$

 \xDH By (C1), repair $ \xbD_{0}:$

$x_{0}= \xCN x_{1} \xcu \xCN x_{2} \xcu \xCN x_{3} \xcu \xCN x_{4},$
$x_{1}= \xCN x_{2} \xcu \xCN x_{3} \xcu \xCN x_{4},$ $x_{2}= \xCN x_{3}
\xcu \xCN x_{4},$ $x_{3}= \xCN x_{4}$

(branch $x_{1} \xcP x_{4}$ has to be contradicted for $x_{0}= \xdT,$ add
$x_{0} \xcP x_{4}).$ Thus

$x_{3}= \xCN x_{4},$ $x_{2}= \xdF,$ $x_{1}$ $=$ $ \xdF,$ $x_{0}$ $=$ $
\xdF.$

$x_{4}$ will be the start of a new Yablo cell.

 \xEj

 \xDH Note:

In the case of triangles $x \xcP y \xcP z,$ $x \xcP z$ formed by paths,
the newly
added $x \xcP z$ must not contain any branchings by ``OR'' after the common
part
of $x \xcP y \xcP z$ and $z \xcP z,$ if $x= \xdT.$ This concerns here
$x_{0} \xcP x_{2}$ for $x_{0}= \xdT $ in (2.2.1),
$x_{1} \xcP x_{3}$ for $x_{1}= \xdT $ in (2.2.2),
$x_{0} \xcP x_{3}$ for $x_{0}= \xdT $ in (2.3),
$x_{2} \xcP x_{4}$ for $x_{2}= \xdT $ in (2.4.2),
$x_{1} \xcP x_{4}$ for $x_{1}= \xdT $ in (2.5),
$x_{0} \xcP x_{4}$ for $x_{0}= \xdT $ in (2.6).
See Remark \ref{Remark Branching} (page \pageref{Remark Branching}), (3)
for explanation and details.

 \xDH Summary (in the limit):

 \xEh
 \xDH Main properties:
 \xEh
 \xDH
Thus, there are arrows $x_{i} \xcP x_{j}$ for all $i,j,$ $i<j,$ and the
construction is
transitive for the $x_{i}' s.$

 \xDH
Every $x_{i}$ is head of a Yablo cell with knee $x_{i+1}.$

 \xDH
Thus, every arrow from any $x_{i}$ to any $x_{j}$ goes to the head of a
Yablo cell, and
not only the arrows from $x_{0}.$

This property is ``accidental'', and due to the fact that for any arrow
$x_{i} \xcP x_{j},$ there is also an arrow $x_{0} \xcP x_{j},$ and (C2)
holds for $x_{0}$
by prerequisite.

 \xEj

 \xDH
The construction has infinite depth and branching.

 \xEj

 \xEj

\ecs

\br

$\hspace{0.01em}$

% {\tiny (+++ Orig. No.:  Remark Arbitrary +++)}

% {\tiny LABEL: {Remark Arbitrary}} \\[1mm]
\label{Remark Arbitrary}

Our construction is somewhat arbitrary. We construct a contradiction for
$x_{0}$ by
$x_{0} \xcP x_{1} \xcP x_{2},$ $x_{0} \xcP x_{2},$ and use the arrow
$x_{1} \xcP x_{2}$ as part of the contradiction for
$x_{1},$ $x_{1} \xcP x_{2} \xcP x_{3},$ $x_{1} \xcP x_{3}.$ This, of
course, is a special case.

We may begin a totally new
contradiction for $x_{1},$ $x_{1} \xcP x_{2}' \xcP x_{3},$ $x_{1} \xcP
x_{3},$ likewise for $x_{2}',$ etc.

But the old construction of
Diagram \ref{Diagram Induction} (page \pageref{Diagram Induction})
is part of this somewhat more general construction,
so the old construction is really minimal.

See also the discussion of ``Saw Blades'' in
Section 7.6 in  \cite{Sch22} and in
Section 3 in  \cite{Sch23b}.

We have to be a bit careful here. If we set $x_{2},$ $x_{3},$ $x_{4},$
etc. $= \xdT,$
then we have $x_{0}= \xdF \xcu \xbx = \xdF,$ etc. Instead, we may set
$x_{2}= \xdF,$ so we have
$x_{0}= \xdT \xcu \xbx = \xbx,$ as desired, or we append at $x_{2}$
recursively the same
construction - as in ``Saw Blades''.

\er

\br

$\hspace{0.01em}$

% {\tiny (+++ Orig. No.:  Remark Reconstruction +++)}

% {\tiny LABEL: {Remark Reconstruction}} \\[1mm]
\label{Remark Reconstruction}

Above construction can be read differently: Suppose we have a
graph $ \xbG $ with valuation such that a node $x_{0}$ has truth value $
\xbx.$
Then we can do above
construction. E.g., $x_{0}= \xdT $ has to be contradictory, (C1), so have
two
contradictory paths $ \xbs,$ $ \xbs' $ from $x_{0}$ to, say $x_{2},$ and
some $x_{1}$ on the path $ \xbs'.$
Both paths have to lead to contradictions for $x_{0}= \xdF,$ (C2).
Exactly one of them
has to be on one of the paths (and not at the end), say this is $x_{1}$
(if
both are on the paths, we would have the diamond).
This way, we can reconstruct the Yablo construction within $ \xbG.$
Of course, there may be additional branchings, as discussed in
Section \ref{Section PathCells} (page \pageref{Section PathCells}).
It will be shown how to neutralise them in
Section \ref{Section Result} (page \pageref{Section Result}).

\clearpage

\begin{diagram}

Diagram Inductive Construction

% {\tiny LABEL: {Diagram Induction}} \\[1mm]
\label{Diagram Induction}
\index{Diagram Induction}

% +BSet

\unitlength0.6mm
\begin{picture}(220,330)(0,0)

% (2.1)

\put(10,310){\circle*{1}}
\put(50,310){\circle*{1}}
\put(90,310){\circle*{1}}
\put(10,310){\line(1,0){80}}

\put(7,305){\footnotesize $x_{0}$}  % (-3,+3)
\put(47,305){\footnotesize $x_{1}$}  % (-3,+3)
\put(87,305){\footnotesize $x_{2}$}  % (-3,+3)
\put(180,305){\footnotesize $(2.1)$}  % (-3,+3)
\put(200,305){\footnotesize $\xbD _0$}  % (-3,+3)

\put(10,302){\line(0,-1){5}}
\put(90,302){\line(0,-1){5}}
\put(10,297){\line(1,0){80}}

\put(49,310){\line(-1,1){1}}
\put(49,310){\line(-1,-1){1}}
\put(89,310){\line(-1,1){1}}
\put(89,310){\line(-1,-1){1}}
\put(90,302){\line(1,-1){1}}
\put(90,302){\line(-1,-1){1}}

\put(30,312){\line(0,-1){4}}
\put(70,312){\line(0,-1){4}}
\put(50,299){\line(0,-1){4}}

% (2.2.1)

\put(10,275){\circle*{1}}
\put(50,275){\circle*{1}}
\put(90,275){\circle*{1}}
\put(130,275){\circle*{1}}
\put(10,275){\line(1,0){120}}

\put(7,270){\footnotesize $x_{0}$}  % (-3,+3)
\put(47,270){\footnotesize $x_{1}$}  % (-3,+3)
\put(87,270){\footnotesize $x_{2}$}  % (-3,+3)
\put(127,270){\footnotesize $x_{3}$}  % (-3,+3)
\put(180,270){\footnotesize $(2.2.1)$, add $\xbD_1$ by (C2) for $x_0 \xcP x_1$}
% (-3,+3)

\put(10,267){\line(0,-1){5}}
\put(90,267){\line(0,-1){5}}
\put(10,262){\line(1,0){80}}

% \put(50,277){\line(0,1){5}}
% \put(130,277){\line(0,1){5}}
% \put(50,282){\line(1,0){80}}

\put(49,275){\line(-1,1){1}}
\put(49,275){\line(-1,-1){1}}
\put(89,275){\line(-1,1){1}}
\put(89,275){\line(-1,-1){1}}
\put(90,267){\line(1,-1){1}}
\put(90,267){\line(-1,-1){1}}

\put(129,275){\line(-1,1){1}}
\put(129,275){\line(-1,-1){1}}

% \put(130,277){\line(-1,1){1}}
% \put(130,277){\line(1,1){1}}

\put(30,277){\line(0,-1){4}}
\put(70,277){\line(0,-1){4}}
\put(50,264){\line(0,-1){4}}
% \put(90,284){\line(0,-1){4}}
\put(110,277){\line(0,-1){4}}

% (2.2.2)

\put(10,240){\circle*{1}}  % -35
\put(50,240){\circle*{1}}
\put(90,240){\circle*{1}}
\put(130,240){\circle*{1}}
\put(10,240){\line(1,0){120}}

\put(7,235){\footnotesize $x_{0}$}  % (-3,+3)
\put(47,235){\footnotesize $x_{1}$}  % (-3,+3)
\put(87,235){\footnotesize $x_{2}$}  % (-3,+3)
\put(127,235){\footnotesize $x_{3}$}  % (-3,+3)
\put(180,235){\footnotesize $(2.2.2)$}  % (-3,+3)

\put(10,232){\line(0,-1){5}}
\put(90,232){\line(0,-1){5}}
\put(10,227){\line(1,0){80}}

\put(50,242){\line(0,1){5}}
\put(130,242){\line(0,1){5}}
\put(50,247){\line(1,0){80}}

\put(49,240){\line(-1,1){1}}
\put(49,240){\line(-1,-1){1}}
\put(89,240){\line(-1,1){1}}
\put(89,240){\line(-1,-1){1}}
\put(90,232){\line(1,-1){1}}
\put(90,232){\line(-1,-1){1}}

\put(130,242){\line(-1,1){1}}
\put(130,242){\line(1,1){1}}

\put(129,240){\line(-1,1){1}}
\put(129,240){\line(-1,-1){1}}

\put(30,242){\line(0,-1){4}}
\put(70,242){\line(0,-1){4}}
\put(110,242){\line(0,-1){4}}
\put(50,229){\line(0,-1){4}}
\put(90,249){\line(0,-1){4}}

% (2.3)

\put(10,200){\circle*{1}}  % -40
\put(50,200){\circle*{1}}
\put(90,200){\circle*{1}}
\put(130,200){\circle*{1}}
\put(10,200){\line(1,0){120}}

\put(7,195){\footnotesize $x_{0}$}  % (-3,+3)
\put(47,195){\footnotesize $x_{1}$}  % (-3,+3)
\put(87,195){\footnotesize $x_{2}$}  % (-3,+3)
\put(127,195){\footnotesize $x_{3}$}  % (-3,+3)
\put(180,195){\footnotesize $(2.3)$, repair $\xbD_0$}
% (-3,+3)

\put(10,192){\line(0,-1){5}}
\put(90,192){\line(0,-1){5}}
\put(10,187){\line(1,0){80}}

\put(50,202){\line(0,1){5}}
\put(128,202){\line(0,1){5}}
\put(50,207){\line(1,0){78}}

\put(10,202){\line(0,1){10}}
\put(132,202){\line(0,1){10}}
\put(10,212){\line(1,0){122}}

\put(49,200){\line(-1,1){1}}
\put(49,200){\line(-1,-1){1}}
\put(89,200){\line(-1,1){1}}
\put(89,200){\line(-1,-1){1}}
\put(90,192){\line(1,-1){1}}
\put(90,192){\line(-1,-1){1}}

\put(128,202){\line(-1,1){1}}
\put(128,202){\line(1,1){1}}

\put(129,200){\line(-1,1){1}}
\put(129,200){\line(-1,-1){1}}

\put(132,202){\line(-1,1){1}}
\put(132,202){\line(1,1){1}}

\put(30,202){\line(0,-1){4}}
\put(70,202){\line(0,-1){4}}
\put(110,202){\line(0,-1){4}}
\put(50,189){\line(0,-1){4}}
\put(90,209){\line(0,-1){4}}
\put(70,214){\line(0,-1){4}}

% (2.4.1)

\put(10,160){\circle*{1}}  % -40
\put(50,160){\circle*{1}}
\put(90,160){\circle*{1}}
\put(130,160){\circle*{1}}
\put(170,160){\circle*{1}}
\put(10,160){\line(1,0){160}}

\put(7,155){\footnotesize $x_{0}$}  % (-3,+3)
\put(47,155){\footnotesize $x_{1}$}  % (-3,+3)
\put(87,155){\footnotesize $x_{2}$}  % (-3,+3)
\put(127,155){\footnotesize $x_{3}$}  % (-3,+3)
\put(167,155){\footnotesize $x_{4}$}  % (-3,+3)
\put(180,155){\footnotesize $(2.4.1)$, add $\xbD_2$ by (C2) for $x_0 \xcP x_2$}
% (-3,+3)

\put(10,152){\line(0,-1){5}}
\put(88,152){\line(0,-1){5}}
\put(10,147){\line(1,0){78}}

\put(50,162){\line(0,1){5}}
\put(128,162){\line(0,1){5}}
\put(50,167){\line(1,0){78}}

\put(10,162){\line(0,1){10}}
\put(132,162){\line(0,1){10}}
\put(10,172){\line(1,0){122}}

% \put(92,152){\line(0,-1){5}}
% \put(170,152){\line(0,-1){5}}
% \put(92,147){\line(1,0){78}}

\put(49,160){\line(-1,1){1}}
\put(49,160){\line(-1,-1){1}}
\put(89,160){\line(-1,1){1}}
\put(89,160){\line(-1,-1){1}}
\put(88,152){\line(1,-1){1}}
\put(88,152){\line(-1,-1){1}}

\put(128,162){\line(-1,1){1}}
\put(128,162){\line(1,1){1}}

\put(129,160){\line(-1,1){1}}
\put(129,160){\line(-1,-1){1}}

\put(132,162){\line(-1,1){1}}
\put(132,162){\line(1,1){1}}

% \put(170,152){\line(1,-1){1}}
% \put(170,152){\line(-1,-1){1}}

\put(30,162){\line(0,-1){4}}
\put(70,162){\line(0,-1){4}}
\put(110,162){\line(0,-1){4}}
\put(50,149){\line(0,-1){4}}
\put(90,169){\line(0,-1){4}}
\put(70,174){\line(0,-1){4}}
% \put(130,149){\line(0,-1){4}}

\put(150,162){\line(0,-1){4}}
\put(169,160){\line(-1,1){1}}
\put(169,160){\line(-1,-1){1}}

% (2.4.2)

\put(10,120){\circle*{1}}  % -40
\put(50,120){\circle*{1}}
\put(90,120){\circle*{1}}
\put(130,120){\circle*{1}}
\put(170,120){\circle*{1}}
\put(10,120){\line(1,0){160}}

\put(7,115){\footnotesize $x_{0}$}  % (-3,+3)
\put(47,115){\footnotesize $x_{1}$}  % (-3,+3)
\put(87,115){\footnotesize $x_{2}$}  % (-3,+3)
\put(127,115){\footnotesize $x_{3}$}  % (-3,+3)
\put(167,115){\footnotesize $x_{4}$}  % (-3,+3)
\put(180,115){\footnotesize $(2.4.2)$}  % (-3,+3)

\put(10,112){\line(0,-1){5}}
\put(88,112){\line(0,-1){5}}
\put(10,107){\line(1,0){78}}

\put(50,122){\line(0,1){5}}
\put(128,122){\line(0,1){5}}
\put(50,127){\line(1,0){78}}

\put(10,122){\line(0,1){10}}
\put(132,122){\line(0,1){10}}
\put(10,132){\line(1,0){122}}

\put(92,112){\line(0,-1){5}}
\put(170,112){\line(0,-1){5}}
\put(92,107){\line(1,0){78}}

\put(49,120){\line(-1,1){1}}
\put(49,120){\line(-1,-1){1}}
\put(89,120){\line(-1,1){1}}
\put(89,120){\line(-1,-1){1}}
\put(88,112){\line(1,-1){1}}
\put(88,112){\line(-1,-1){1}}

\put(128,122){\line(-1,1){1}}
\put(128,122){\line(1,1){1}}

\put(129,120){\line(-1,1){1}}
\put(129,120){\line(-1,-1){1}}

\put(132,122){\line(-1,1){1}}
\put(132,122){\line(1,1){1}}

\put(170,112){\line(1,-1){1}}
\put(170,112){\line(-1,-1){1}}

\put(169,120){\line(-1,1){1}}
\put(169,120){\line(-1,-1){1}}

\put(30,122){\line(0,-1){4}}
\put(70,122){\line(0,-1){4}}
\put(110,122){\line(0,-1){4}}
\put(150,122){\line(0,-1){4}}
\put(50,109){\line(0,-1){4}}
\put(90,129){\line(0,-1){4}}
\put(70,134){\line(0,-1){4}}
\put(130,109){\line(0,-1){4}}

% (2.5)

\put(10,80){\circle*{1}}  % -40
\put(50,80){\circle*{1}}
\put(90,80){\circle*{1}}
\put(130,80){\circle*{1}}
\put(170,80){\circle*{1}}
\put(10,80){\line(1,0){160}}

\put(7,75){\footnotesize $x_{0}$}  % (-3,+3)
\put(47,75){\footnotesize $x_{1}$}  % (-3,+3)
\put(87,75){\footnotesize $x_{2}$}  % (-3,+3)
\put(127,75){\footnotesize $x_{3}$}  % (-3,+3)
\put(167,75){\footnotesize $x_{4}$}  % (-3,+3)
\put(180,75){\footnotesize $(2.5)$, repair $\xbD_1$, but damage $\xbD_0$}
% (-3,+3)

\put(10,72){\line(0,-1){5}}
\put(88,72){\line(0,-1){5}}
\put(10,67){\line(1,0){78}}

\put(50,82){\line(0,1){5}}
\put(128,82){\line(0,1){5}}
\put(50,87){\line(1,0){78}}

\put(10,82){\line(0,1){10}}
\put(132,82){\line(0,1){10}}
\put(10,92){\line(1,0){122}}

\put(92,72){\line(0,-1){5}}
\put(166,72){\line(0,-1){5}}
\put(92,67){\line(1,0){74}}

\put(50,72){\line(0,-1){10}}
\put(170,72){\line(0,-1){10}}
\put(50,62){\line(1,0){120}}

\put(49,80){\line(-1,1){1}}
\put(49,80){\line(-1,-1){1}}
\put(89,80){\line(-1,1){1}}
\put(89,80){\line(-1,-1){1}}
\put(88,72){\line(1,-1){1}}
\put(88,72){\line(-1,-1){1}}

\put(128,82){\line(-1,1){1}}
\put(128,82){\line(1,1){1}}

\put(129,80){\line(-1,1){1}}
\put(129,80){\line(-1,-1){1}}

\put(132,82){\line(-1,1){1}}
\put(132,82){\line(1,1){1}}

\put(166,72){\line(1,-1){1}}
\put(166,72){\line(-1,-1){1}}

\put(170,72){\line(1,-1){1}}
\put(170,72){\line(-1,-1){1}}

\put(169,80){\line(-1,1){1}}
\put(169,80){\line(-1,-1){1}}

\put(30,82){\line(0,-1){4}}
\put(70,82){\line(0,-1){4}}
\put(110,82){\line(0,-1){4}}
\put(150,82){\line(0,-1){4}}
\put(45,69){\line(0,-1){4}}
\put(90,89){\line(0,-1){4}}
\put(70,94){\line(0,-1){4}}
\put(130,69){\line(0,-1){4}}
\put(110,64){\line(0,-1){4}}

% (2.6)

\put(10,36){\circle*{1}}  % -40
\put(50,36){\circle*{1}}
\put(90,36){\circle*{1}}
\put(130,36){\circle*{1}}
\put(170,36){\circle*{1}}
\put(10,36){\line(1,0){160}}

\put(7,31){\footnotesize $x_{0}$}  % (-3,+3)
\put(47,31){\footnotesize $x_{1}$}  % (-3,+3)
\put(87,31){\footnotesize $x_{2}$}  % (-3,+3)
\put(127,31){\footnotesize $x_{3}$}  % (-3,+3)
\put(167,31){\footnotesize $x_{4}$}  % (-3,+3)
\put(180,31){\footnotesize $(2.6)$, repair $\xbD_0$}
% (-3,+3)

\put(10,28){\line(0,-1){5}}
\put(88,28){\line(0,-1){5}}
\put(10,23){\line(1,0){78}}

\put(50,38){\line(0,1){5}}
\put(128,38){\line(0,1){5}}
\put(50,43){\line(1,0){78}}

\put(10,38){\line(0,1){10}}
\put(132,38){\line(0,1){10}}
\put(10,48){\line(1,0){122}}

\put(92,28){\line(0,-1){5}}
\put(166,28){\line(0,-1){5}}
\put(92,23){\line(1,0){74}}

\put(50,28){\line(0,-1){10}}
\put(170,28){\line(0,-1){10}}
\put(50,18){\line(1,0){120}}

\put(8,28){\line(0,-1){15}}
\put(174,28){\line(0,-1){15}}
\put(8,13){\line(1,0){166}}

\put(49,36){\line(-1,1){1}}
\put(49,36){\line(-1,-1){1}}
\put(89,36){\line(-1,1){1}}
\put(89,36){\line(-1,-1){1}}
\put(88,28){\line(1,-1){1}}
\put(88,28){\line(-1,-1){1}}

\put(128,38){\line(-1,1){1}}
\put(128,38){\line(1,1){1}}

\put(129,36){\line(-1,1){1}}
\put(129,36){\line(-1,-1){1}}

\put(132,38){\line(-1,1){1}}
\put(132,38){\line(1,1){1}}

\put(166,28){\line(1,-1){1}}
\put(166,28){\line(-1,-1){1}}

\put(170,28){\line(1,-1){1}}
\put(170,28){\line(-1,-1){1}}

\put(169,36){\line(-1,1){1}}
\put(169,36){\line(-1,-1){1}}

\put(174,28){\line(-1,-1){1}}
\put(174,28){\line(1,-1){1}}

\put(30,38){\line(0,-1){4}}
\put(70,38){\line(0,-1){4}}
\put(110,38){\line(0,-1){4}}
\put(150,38){\line(0,-1){4}}
\put(45,25){\line(0,-1){4}}
\put(90,45){\line(0,-1){4}}
\put(70,50){\line(0,-1){4}}
\put(130,25){\line(0,-1){4}}
\put(110,20){\line(0,-1){4}}
\put(90,15){\line(0,-1){4}}

\end{picture}

\end{diagram}

\vspace{10mm}

\clearpage

\clearpage
\subsection{
Variant Saw-Blade
}

% {\tiny LABEL: {Section Combining-Saw-Blade}} \\[1mm]
\label{Section Combining-Saw-Blade}
\subsubsection{
Saw Blades
}

% {\tiny LABEL: {Section Saw-Blades-All}} \\[1mm]
\label{Section Saw-Blades-All}

\er

For the convenience of the reader, we first recall the main elements of
the Saw
Blade construction, as published in
 \cite{Sch21}, v3, or  \cite{Sch22}.

Our present construction,
see Section 
\ref{Section Induction-New-1} (page 
\pageref{Section Induction-New-1})
discusses one of the simplified versions, but in much more detail,
shows that continuation via the head does not suffice, see (5) there,
and shows the necessity of the details of the construction in case
(1.2.2) of
Section \ref{Section Combining-Intro} (page \pageref{Section Combining-Intro}).

\bcs

$\hspace{0.01em}$

% {\tiny (+++ Orig. No.:  Construction Saw-Blade +++)}

% {\tiny LABEL: {Construction Saw-Blade}} \\[1mm]
\label{Construction Saw-Blade}

We construct a saw blade $ \xbs,$ $SB_{ \xbs }.$

 \xEh

 \xDH
``Saw Blades''
 \xEh
 \xDH
Let $x_{ \xbs,0} \xcP x_{ \xbs,1} \xcP x_{ \xbs,2} \xcP x_{ \xbs,3}
\xcP x_{ \xbs,4},$  \Xl.

$x_{ \xbs,0} \xcP y_{ \xbs,0},$ $x_{ \xbs,1} \xcP y_{ \xbs,0},$ $x_{
\xbs,1} \xcP y_{ \xbs,1},$ $x_{ \xbs,2} \xcP y_{ \xbs,1},$
$x_{ \xbs,2} \xcP y_{ \xbs,2},$ $x_{ \xbs,3} \xcP y_{ \xbs,2},$ $x_{
\xbs,3} \xcP y_{ \xbs,3},$ $x_{ \xbs,4} \xcP y_{ \xbs,3},$  \Xl.

we call the construction
a ``saw blade'', with ``teeth'' $y_{ \xbs,0},$ $y_{ \xbs,1},$ $y_{ \xbs
,2},$  \Xl.
and ``back'' $x_{ \xbs,0},$ $x_{ \xbs,1},$ $x_{ \xbs,2},$  \Xl.

We call $x_{ \xbs,0}$ the start of the blade.

 \xDH
Add (against escape), e.g. first $x_{ \xbs,0} \xcP x_{ \xbs,2},$ $x_{
\xbs,0} \xcP y_{ \xbs,1},$
then $x_{ \xbs,1} \xcP x_{ \xbs,3},$ $x_{ \xbs,1} \xcP y_{ \xbs,2},$
now
we have to add $x_{ \xbs,0} \xcP x_{ \xbs,3},$ $x_{ \xbs,0} \xcP y_{
\xbs,2},$ etc, recursively.
This is equivalent to
closing the saw blade under transitivity with negative arrows $ \xcP.$
This is easily seen.

 \xDH
We set $x_{ \xbs,i}= \xcU \xCN z_{ \xbs,j},$ for all $x_{ \xbs,j}$
such that $x_{ \xbs,i} \xcP z_{ \xbs,j},$
as in the original Yablo construction.

 \xEj

 \xDH
Composition of saw blades
 \xEh
 \xDH
Add for the teeth $y_{ \xbs,0},$ $y_{ \xbs,1},$ $y_{ \xbs,2}$  \Xl.
their own saw blades, i.e.
start at $y_{ \xbs,0}$ a new saw blade $SB_{ \xbs,0}$ with $y_{ \xbs
,0}=x_{ \xbs,0,0},$ at
$y_{ \xbs,1}$ a new saw blade $SB_{ \xbs,1}$ with $y_{ \xbs,1}=x_{ \xbs
,1,0},$ etc.

 \xDH
Do this recursively.

I.e., at every tooth of every saw blade start a new saw blade.
 \xEj
 \xEj

Note:

It is NOT necessary to close the whole structure (the individual saw
blades
together) under transitivity.
\paragraph{
Simplifications of the Saw Blade Construction
}

% {\tiny LABEL: {Section Simplifications}} \\[1mm]
\label{Section Simplifications}

\ecs

We show here that it is not necessary to make the $y_{ \xbs,i}$
contradictory
in a recursive construction, as in
Construction 
\ref{Construction Saw-Blade} (page 
\pageref{Construction Saw-Blade}).
It suffices to prevent them to be true.

We discuss three, much simplified, Saw Blade constructions.

Note, however, that the back of each saw blade ``hides'' a Yablo
construction. The separate treatment of the teeth illustrates the
conceptual difference, but it cannot escape blurring it again in
the back of the blade.

\bcs

$\hspace{0.01em}$

% {\tiny (+++ Orig. No.:  Construction Simple-1 +++)}

% {\tiny LABEL: {Construction Simple-1}} \\[1mm]
\label{Construction Simple-1}

 \xEh

 \xDH

Take ONE saw blade $ \xbs,$ and attach (after closing under transitivity)
at all $y_{ \xbs,i}$ a SINGLE Yablo Cell $y_{ \xbs,i} \xcP u_{ \xbs,i}
\xcP v_{ \xbs,i},$ $y_{ \xbs,i} \xcP v_{ \xbs,i}.$
We call this the decoration, it is not involved in closure under
transitivity.

 \xEh

 \xDH

Any node $z$ in the saw blade (back or tooth) cannot be $z+,$ this leads
to a contradiction:

If $z=x_{i}$ (in the back):

Let $x_{i}+:$
Take $x_{i} \xcP x_{i+1}$ (any $x_{j},$ $i<j$ would do), if $x_{i+1} \xcP
r,$ then by transitivity,
$x_{i} \xcP r,$ so we have a contradiction.

If $z=y_{i}$ (a tooth):

$y_{i}+$ is contradictory by the ``decoration'' appended to $y_{i}.$

 \xDH

Any $x_{i}-$ $(x_{i}$ in the back, as a matter of fact, $x_{0}-$ would
suffice)
is impossible:

Consider any $x_{i} \xcP r,$ then $r+$ is impossible, as we just saw.

Note: there are no arrows from the back of the blade to the decoration.

 \xEj

 \xDH

We can simplify even further. The only thing we need about the $y_{i}$ is
that
they cannot be $+.$ Instead of decorating them with a Yablo Cell, any
contradiction will do, the simplest one is $y_{i}=y'_{i} \xcu \xCN
y'_{i}.$ Even
just one $y' $ s.t. $y_{i}=y' \xcu \xCN y' $ for all $y_{i}$ would do. (Or
a constant
FALSE.)

Formally, we set

$x_{0}:= \xcU \{ \xCN x_{i}:i>0\} \xcu \xcU \{ \xCN y_{i}:i \xcg 0\},$

for $j>0:$

$x_{j}:= \xcU \{ \xCN x_{i}:i>j\} \xcu \xcU \{ \xCN y_{i}:i \xcg j-1\},$

and

$y_{j}:=y'_{j} \xcu \xCN y'_{j}.$

 \xDH

In a further step, we see that the $y_{i}$ (and thus the $y'_{i})$ need
not
be different from each other, one $y$ and one $y' $ suffice.

Thus, we set $x_{j}:= \xcU \{ \xCN x_{i}:i>j\} \xcu \xCN y,$ $y:=y' \xcu
\xCN y'.$

(Intuitively, the cells are arranged in a circle, with $y$ at the center,
and
$y' $ ``sticking out''. We might call this a ``curled saw blade''.

 \xDH
When we throw away the $y_{j}$ altogether, we have Yablo's construction.
this works, as we have the essential part in the $x_{i}$'s, and used the
$y_{j}' s$
only as a sort of scaffolding.

 \xEj
\subsubsection{
Our present construction and its discussion
}

% {\tiny LABEL: {Section Induction-New-1}} \\[1mm]
\label{Section Induction-New-1}

\ecs

\br

$\hspace{0.01em}$

% {\tiny (+++ Orig. No.:  Remark Induction-New-1 +++)}

% {\tiny LABEL: {Remark Induction-New-1}} \\[1mm]
\label{Remark Induction-New-1}

 \xEh

 \xDH The construction is very similar to the construction in
Section 
\ref{Section Combining-Yablo} (page 
\pageref{Section Combining-Yablo}), just
a bit more complicated.

 \xDH In Diagram 
\ref{Diagram Induction-New} (page 
\pageref{Diagram Induction-New}),
we omit the added arrows like $x_{0} \xcP x_{2}$ for better readability.
(Such arrows were shown in
Diagram \ref{Diagram Induction} (page \pageref{Diagram Induction}).)

 \xDH We note $a_{i,j}:x_{i} \xcP x_{j},$ $b_{i,j}:x_{i} \xcP y_{j},$ and
$c_{i,j}$ the arrows in the Upper Part
of Diagram \ref{Diagram Induction-New} (page \pageref{Diagram Induction-New}).

 \xDH Basically, we need (C2) of $x_{0}$ to extend the construction,
and preservation of contradiction in $ \xbD_{i}$ (see below) for
transitivity - just as in
Diagram \ref{Diagram Induction} (page \pageref{Diagram Induction}).

 \xEh

 \xDH Let $ \xbD_{i}$ $:=$ $x_{i} \xcP x_{i+1} \xcP y_{i},$ $x_{i} \xcP
y_{i},$ we abbreviate $x_{0} \xcP_{ \xbD_{0}}x_{1},$ etc.

 \xDH $x_{0} \xcP_{ \xbD_{0}}x_{1} \xcP_{ \xbD_{1}}x_{2} \xcP_{
\xbD_{2}}x_{3} \xcP_{ \xbD_{3}}x_{4}.$

 \xDH $ \xbD_{0}$ has to exist by (C1) fuer $x_{0}.$

 \xDH $ \xbD_{1}$ has to exist by (C2) fuer $x_{0}.$

 \xDH $a_{0,2}:x_{0} \xcP x_{2}$ repairs $ \xbD_{0},$ restoring the
contradiction, so this has to exist.
For the same reason, there is $b_{0,1' }:x_{0} \xcP y_{1}.$
From now on, we will not mention the $b_{i,j}$ any more, for simplicity's
sake.

 \xDH By (C2) for $x_{0}$ and $a_{0,2},$ we need $ \xbD_{2}.$

 \xDH $a_{1,3}:x_{1} \xcP x_{3}$ repairs $ \xbD_{1},$ and damages $
\xbD_{0},$ $a_{0,3}$ is needed to repair $ \xbD_{0}.$

These constructions result in transitivity along the x's line, so we
see that the construction contains the Yablo construction.

 \xDH By (C2) for $x_{0}$ and $a_{0,3}$ we need $ \xbD_{3},$ etc.

 \xEj

 \xDH We now show that attaching new contradiction cells at the head,
and not at the knee of the old cell does not suffice. Basically,
as the old contradiction is not finished at the knee, we have to
continue, at the head, the old cell is finished, and there is no
need to connect the old cells to the new cell, thus interrupting
the construction. This shows the importance of the distinction
between head and knee.

 \xEh

 \xDH Take $x_{0} \xcP_{ \xbD_{0}}x_{1} \xcP_{ \xbD_{1}}x_{2} \xcP_{
\xbD_{2}}x_{3} \xcP_{ \xbD_{3}}x_{4}.$

Add

$y_{2} \xcP z_{1} \xcP z_{2},$ $y_{2} \xcP z_{2}$ und $z_{1} \xcP z_{1,1}
\xcP z_{1,2},$ the Upper Part of
Diagram \ref{Diagram Induction-New} (page \pageref{Diagram Induction-New}).

 \xDH The difference: $x_{3}$ is a knee, to preserve $ \xbD_{2}$
contradictory, we need $a_{2,4}$
to pass over $x_{3}$ hinweg. This is not the case with $y_{2},$ as it is a
head, not
a knee.

 \xDH Thus, without $x_{4},$ $y_{3}$ the construction ends at $y_{2},$ the
lower part is
finished.

More precisely:

 \xEh
 \xDH (C2) for $x_{0}$ is satisfied at $y_{2},$ so is (C1) for $
\xbD_{2},$ da $y_{2}$ is a head.
If $x_{4},$ $y_{3}$ do not exist, $ \xbD_{2}$ is also satisfied at
$x_{3}.$
 \xDH $a_{0,2}:x_{0} \xcP x_{2},$ $a_{0,3}:x_{0} \xcP x_{3},$
$a_{1,3}:x_{1} \xcP x_{3},$ and
$b_{0,1}:x_{0} \xcP y_{1},$ $b_{0,2}:x_{0} \xcP y_{2},$ $b_{1,2}:x_{1}
\xcP y_{2}$
satify all (C1) conditions, so this is closed.
 \xDH We do not need any arrows from the lower part to the upper part
 \xDH Finally, we may calculate $x_{3},$ $x_{2},$ $x_{1},$ $x_{0}$ from
$y_{0},$ $y_{1},$ $y_{2}.$
So this is only a classical addition.
 \xEj

 \xDH Summary
 \xEh
 \xDH
The boundary is $y_{2},$ it separates the lower part (below $y_{2})$ from
the
upper part, we are not forced to cross from the lower to the upper part.

 \xDH
If we append at $x_{3}$ a FINITE construction, like $ \xdF,$ to satisfy
(C2) and
(C2), we may fill in from below some values, and calculate the beginning
of the diagram, resulting in classical values.
So, the start of the diagram is then isolated from the rest, and thus
classical.
 \xEj

 \xEj

 \xDH Suppose $ \xbD_{2}$ is a YC until $x_{3},$ but we do not defend it
against
(the necessary) $ \xbD_{3}.$

 \xEh

 \xDH $x_{0} \xcP_{ \xbD_{0}}x_{1} \xcP_{ \xbD_{1}}x_{2} \xcP_{
\xbD_{2}}x_{3} \xcP_{ \xbD_{3}}x_{4}.$

 \xDH Then there might not be some $a_{2,4}:x_{2} \xcP x_{4},$ so the
possibility at
$x_{3},$ $x_{3} \xcP x_{4}$ is not contradicted (by $x_{2} \xcP x_{4}),$
and $ \xbD_{2}$ is destroyed.

 \xDH It does not help to consider a group of incomplete $ \xbD_{2}' s,$
in any case the $x_{3,i}$ keep a destructive possibility, which may
be chosen by the arbitrary choice of the $x_{4,i}.$ We need a complete $
\xbD_{2},$
which stays contradictory.

 \xEj

 \xEj

\er

\br

$\hspace{0.01em}$

% {\tiny (+++ Orig. No.:  Remark Induction-New-2 +++)}

% {\tiny LABEL: {Remark Induction-New-2}} \\[1mm]
\label{Remark Induction-New-2}

The details are very similar to the case discussed in
Section 
\ref{Section Combining-Yablo} (page 
\pageref{Section Combining-Yablo}), so we can be concise.

Consider
Diagram \ref{Diagram Induction-New} (page \pageref{Diagram Induction-New})
(without the upper part, which only served to illustrate the importance of
the construction via the knees).

The ``repair steps'' are needed to repair contradictions damaged
by newly added possibilities.
(See the end of
Remark \ref{Remark Induction-New-1} (page \pageref{Remark Induction-New-1}).)

 \xEh
 \xDH Start with $ \xbD_{0},$ (C1) for $x_{0}.$
 \xDH Add $ \xbD_{1}$ by (C2) and $x_{0} \xcP x_{1},$ damage $ \xbD_{0}.$
 \xDH Repair $ \xbD_{0}$ with $a_{0,2}.$
 \xDH Add $ \xbD_{2}$ by (C2) and $a_{0,2},$ damage $ \xbD_{1}.$
 \xDH Repair $ \xbD_{1}$ with $a_{1,3},$ damage $ \xbD_{0}.$
 \xDH Repair $ \xbD_{0}$ with $a_{0,3}.$
 \xDH Add $ \xbD_{3}$ by (C2) and $a_{0,3},$ damage $ \xbD_{2}.$
 \xDH Repair $ \xbD_{2}$ with $a_{2,4},$ damage $ \xbD_{1}.$
 \xDH Repair $ \xbD_{1}$ with $a_{1,4},$ damage $ \xbD_{0}.$
 \xDH Repair $ \xbD_{0}$ with $a_{0,4}.$
 \xDH Add $ \xbD_{4},$ etc.
 \xEj

\clearpage
\begin{diagram}

% {\tiny LABEL: {Diagram Induction-New}} \\[1mm]
\label{Diagram Induction-New}
\index{Diagram New}

\unitlength.6mm
\begin{picture}(280,320)(0,0)

% x

\put(40,40){\circle*{1}}
\put(80,80){\circle*{1}}
\put(120,120){\circle*{1}}
\put(160,160){\circle*{1}}
\put(200,200){\circle*{1}}
\put(240,240){\circle*{1}}

\put(38,30){\footnotesize $x_0$}
\put(85,75){\footnotesize $x_1$}
\put(125,115){\footnotesize $x_2$}
\put(165,155){\footnotesize $x_3$}
\put(205,195){\footnotesize $x_4$}
\put(245,235){\footnotesize $x_5$}

% \put(45,30){\footnotesize $40,40$}
% \put(92,75){\footnotesize $80,80$}
% \put(132,115){\footnotesize $120,120$}
% \put(172,155){\footnotesize $160,160$}
% \put(212,195){\footnotesize $200,200$}
% \put(238,248){\footnotesize $240,240$}

\put(40,40){\line(1,1){120}}
\put(200,200){\line(1,1){40}}

\put(80,80){\line(-1,0){2}}
\put(80,80){\line(0,-1){2}}
\put(120,120){\line(-1,0){2}}
\put(120,120){\line(0,-1){2}}
\put(160,160){\line(-1,0){2}}
\put(160,160){\line(0,-1){2}}
\put(160,160){\line(1,1){40}}
\put(200,200){\line(-1,0){2}}
\put(200,200){\line(0,-1){2}}
\put(240,240){\line(-1,0){2}}
\put(240,240){\line(0,-1){2}}

%y

\put(40,120){\circle*{1}}
\put(80,160){\circle*{1}}
\put(120,200){\circle*{1}}
\put(160,240){\circle*{1}}
\put(200,280){\circle*{1}}

\put(35,125){\footnotesize $y_0$}
\put(75,165){\footnotesize $y_1$}
\put(110,205){\footnotesize $y_2$}
\put(155,245){\footnotesize $y_3$}
\put(195,285){\footnotesize $y_4$}

% \put(42,125){\footnotesize $40,120$}
% \put(82,165){\footnotesize $80,160$}
% \put(117,205){\footnotesize $120,200$}
% \put(162,245){\footnotesize $160,240$}
% \put(202,285){\footnotesize $200,280$}
%

\put(40,40){\line(0,1){80}}
\put(80,80){\line(0,1){80}}
\put(120,120){\line(0,1){80}}
\put(160,160){\line(0,1){80}}
\put(200,200){\line(0,1){80}}

\put(40,118){\line(-1,-1){2}}
\put(40,118){\line(1,-1){2}}
\put(80,158){\line(-1,-1){2}}
\put(80,158){\line(1,-1){2}}
\put(120,198){\line(-1,-1){2}}
\put(120,198){\line(1,-1){2}}
\put(160,238){\line(-1,-1){2}}
\put(160,238){\line(1,-1){2}}
\put(200,278){\line(-1,-1){2}}
\put(200,278){\line(1,-1){2}}

\put(43,117){\line(0,-1){2}}
\put(43,117){\line(1,0){2}}
\put(83,157){\line(0,-1){2}}
\put(83,157){\line(1,0){2}}
\put(123,197){\line(0,-1){2}}
\put(123,197){\line(1,0){2}}
\put(163,237){\line(0,-1){2}}
\put(163,237){\line(1,0){2}}
\put(203,277){\line(0,-1){2}}
\put(203,277){\line(1,0){2}}

\put(80,80){\line(-1,1){40}}
\put(120,120){\line(-1,1){40}}
\put(160,160){\line(-1,1){40}}
\put(200,200){\line(-1,1){40}}
\put(240,240){\line(-1,1){40}}

% Sonderteil

\put(120,230){\circle*{1}}
\put(120,270){\circle*{1}}
\put(140,250){\circle*{1}}
\put(100,250){\circle*{1}}

\put(125,228){\footnotesize $z_{1}$}
\put(90,250){\footnotesize $z_{2}$}
\put(142,250){\footnotesize $z_{1,1}$}
\put(115,275){\footnotesize $z_{1,2}$}

\put(120,200){\line(0,1){70}}
\put(120,230){\line(-1,1){20}}
\put(120,230){\line(1,1){20}}
\put(140,250){\line(-1,1){20}}
\put(120,200){\line(-2,5){20}}

\put(120,230){\line(-1,-1){2}}
\put(120,230){\line(1,-1){2}}
\put(120,268){\line(-1,-1){2}}
\put(120,268){\line(1,-1){2}}

\put(140,250){\line(0,-1){2}}
\put(140,250){\line(-1,0){2}}
\put(103,247){\line(0,-1){2}}
\put(103,247){\line(1,0){2}}
\put(123,267){\line(0,-1){2}}
\put(123,267){\line(1,0){2}}

\put(102,245){\line(0,-1){3}}
\put(102,245){\line(1,-1){2}}

\put(55,80){\footnotesize $\xbD_0$}
\put(95,120){\footnotesize $\xbD_1$}
\put(135,160){\footnotesize $\xbD_2$}
\put(175,200){\footnotesize $\xbD_3$}
\put(215,240){\footnotesize $\xbD_4$}

\put(111,225){\footnotesize $\xbD_{y_2}$}
\put(125,249){\footnotesize $\xbD_{z_{1}}$}

\multiput(152,202)(-2,0){35}{\circle*{0.3}}
\multiput(152,202)(0,2){40}{\circle*{0.3}}

\put(80,220){\footnotesize Upper}
\put(80,215){\footnotesize Part}

\end{picture}

\end{diagram}

\vspace{10mm}

\clearpage

\clearpage
\clearpage
\subsection{
Conclusion
}

% {\tiny LABEL: {Section Con-3}} \\[1mm]
\label{Section Con-3}

\er

We have shown here (together with the results on Saw Blades,
see Section 7.6 in  \cite{Sch22} and
Section 3 in  \cite{Sch23b})
that - within our framework -
the universal way to combine contradictions to have $ \xfA x \xfA = \xbx $
for some $x$
is Yablo's construction. Other constructions (Saw Blades) are possible,
but they contain Yablo's Construction.

For this reason, whether constructed ``from scratch'', or within an
existing
structure $ \xbG $ with some $x$ with $ \xfA x \xfA = \xbx,$ we find
Yablo's structure.

This is essential for our final result.
\clearpage
\section{
Result
}

% {\tiny LABEL: {Section Result}} \\[1mm]
\label{Section Result}

We use here our results in
Section 
\ref{Section Elementary-Cells} (page 
\pageref{Section Elementary-Cells})  and
Section \ref{Section Combining-Cells} (page \pageref{Section Combining-Cells}),
i.e. the universality
of Yablo's construction, to show that any acyclic directed graph with
an element $x$ such that both $x$ and $ \xCN x$ lead to a contradiction
must
contain a suitable image of Yablo's construction, and conversely,
the existence of such an image guarantees the existence of such an
element.

\br

$\hspace{0.01em}$

% {\tiny (+++ Orig. No.:  Remark Neutral +++)}

% {\tiny LABEL: {Remark Neutral}} \\[1mm]
\label{Remark Neutral}

We use here the following (trivial) idea:

Let $ \xbG =\{x \xcp x',$ $x \xcp x'' \}$ with valuation (for $x)$ e.g. $
\xbf_{x}=x' \xcu \xCN x''.$
We add now $x \xcp y$ to $ \xbG,$ and want to adjust $ \xbf_{x}$ in a
neutral way.
We have to add something by
Definition \ref{Definition Graph} (page \pageref{Definition Graph}), (3.1).
Adding $y \xco \xCN y$ with the result $(x' \xcu \xCN x'') \xcu (y \xco
\xCN y)$
takes care of this, but, if $ \xfA y \xfA = \xbx,$ this is not neutral,
so
we add $y \xco \xCN y \xco \xdT.$
$x_{i}$ $=$ $x_{i} \xcu (x_{j} \xco \xCN x_{j} \xco \xdT),$ or its
converse $x_{i}$ $=$ $x_{i} \xco (x_{j} \xcu \xCN x_{j} \xcu \xdF).$

More generally, if $ \xbf_{x}$ is given for $ \xbG,$ and we add some $x
\xcp x_{j},$ $j \xbe J,$ to $ \xbG,$ then
we consider $ \xbf_{x} \xcu \xcU \{x_{j} \xco \xCN x_{j} \xco \xdT:$ $j
\xbe J\}.$

By its logical neutrality, the addition does not interfere with the
constructions.

\er

\paragraph{
The Problem
}

$\hspace{0.01em}$

% {\tiny (+++ Orig.:  The Problem +++)}

% {\tiny LABEL: {Section The Problem}} \\[1mm]
\label{Section The Problem}

Given an acyclic directed graph $ \xbG,$ does it have a valuation of the
syntactic form $ \xcU \xCN x_{i}$ or $ \xcU \xCN x_{i} \xcu (y \xco \xCN y
\xco \xdT)$ such that
there is a node with truth value $ \xbx?$

We first show that, if there is a suitable injection from the Yablo
structure
YS into $ \xbG,$ then $ \xbG $ has a valuation as above with some $x,$ $
\xfA x \xfA = \xbx,$ see
Section \ref{Section Injection} (page \pageref{Section Injection}).

We then show, if a valuation as above exists for $ \xbG,$ such that for
some
$x$ $ \xfA x \xfA = \xbx,$ then we can reconstruct YS inside $ \xbG,$
and thus have again
a suitable injection, see
Section \ref{Section Converse} (page \pageref{Section Converse}).

We did the main work
in Section 
\ref{Section Combining-Cells} (page 
\pageref{Section Combining-Cells}).
We showed that all constructions with the Yablo Property
(see Definition 
\ref{Definition Yablo-Structure} (page 
\pageref{Definition Yablo-Structure}))
contain the Yablo structure, Yablo has shown that this
structure has the Yablo Property.

So the question is: does a graph contain (a suitable image) of the
Yablo Structure? This is now easy.
\subsection{
A Suitable Injection From the Yablo Structure Proves the Existence
of some $x$ such that $x$ and $\xCN x$ are Contradictory
}

% {\tiny LABEL: {Section Injection}} \\[1mm]
\label{Section Injection}

Consider Yablo's original structure $ \xCf YS,$ and some directed acyclic
graph $ \xbG $
without valuation.

Then there is a valuation of $ \xbG $ in the extended language as
described
in Remark \ref{Remark Neutral} (page \pageref{Remark Neutral})
which gives some element $x_{0}$ the truth value $ \xbx $
if the following conditions hold.

\bcd

$\hspace{0.01em}$

% {\tiny (+++ Orig. No.:  Condition Injection +++)}

% {\tiny LABEL: {Condition Injection}} \\[1mm]
\label{Condition Injection}

 \xEh
 \xDH There is an injection $ \xbm $ from $ \xCf YS$ to $ \xbG,$ more
precisely from
the nodes and paths of $ \xCf YS$ to the nodes and paths of $ \xbG,$
s.t.,
if there is a path $ \xbs $ from $x$ to $y$ in $ \xCf YS,$
then there is a path $ \xbm (\xbs)$ from $ \xbm (x)$ to $ \xbm (y)$ in $
\xbG,$
and if $x' $ is an intermediate node on $ \xbs,$ $ \xbm (x')$ will be an
intermediate
node on $ \xbm (\xbs)$ (but $ \xbm (\xbs)$ may have additional
intermediate nodes).

 \xDH $ \xbm $ preserves $ \xCf pol$ on paths: $pol(\xbs)=pol(\xbm (
\xbs)).$
See Definition 
\ref{Definition Path-Value} (page 
\pageref{Definition Path-Value}).

Thus, if $ \xbs $ contradicts $ \xbs',$ so does $ \xbm (\xbs)$ to $
\xbm (\xbs').$

The second condition expresses ``richness'',
see e.g.  \cite{Sch23b},
that we have enough
independent paths.

 \xEj

\ecd

\be

$\hspace{0.01em}$

% {\tiny (+++ Orig. No.:  Example Pol +++)}

% {\tiny LABEL: {Example Pol}} \\[1mm]
\label{Example Pol}

Let $ \xbs =x_{0} \xcP x_{1},$ $ \xbs' =x_{1} \xcP x_{2},$ $ \xbt =x_{0}
\xcP x_{2},$
$ \xbm (\xbt)= \xbm (x_{0}) \xcP x'_{1} \xcP \xbm (x_{2}),$ then the
second condition fails for $ \xbt.$
(It would not fail for $ \xbm (\xbt)= \xbm (x_{0}) \xcP x'_{1} \xcP
x''_{1} \xcP \xbm (x_{2}).)$

\ee

Let $ \xbm [YS]$ be the image of YS in $ \xbG.$
Due to the syntactic restriction of the $ \xbf_{x}$'s, $ \xbm [YS]$ has a
unique valuation,
which is a partial valuation for $ \xbG.$

\bd

$\hspace{0.01em}$

% {\tiny (+++ Orig. No.:  Definition Extension +++)}

% {\tiny LABEL: {Definition Extension}} \\[1mm]
\label{Definition Extension}

We extend this partial valuation in $ \xbm [YS]$ to a total valuation in $
\xbG $ as
follows:

 \xEh
 \xDH If $x$ is a node in $ \xbG - \xbm [YS],$ set
$ \xbf_{x}$ $=$ $ \xcU \{(y \xco \xCN y \xco \xdT):$ there is an arrow $x
\xcp y$ in $ \xbG \}.$
 \xDH If $x$ is a node in $ \xbm [YS]:$
We have already a formula $ \xbf_{x}' $ describing the behaviour for all
successors of $x$ in $ \xbm [YS].$ We now have to respect the other arrows
in $ \xbG,$
and set

$ \xbf_{x}$ $=$ $ \xbf_{x}' $ $ \xcu $ $ \xcU \{(y \xco \xCN y \xco \xdT
):$ there is an arrow $x \xcp y$ in $ \xbG,$ and $y$
not a successor of $x$ in $ \xbm [YS]$ $\}.$
 \xEj

\ed

\bfa

$\hspace{0.01em}$

% {\tiny (+++ Orig. No.:  Fact Injection +++)}

% {\tiny LABEL: {Fact Injection}} \\[1mm]
\label{Fact Injection}

By Condition \ref{Condition Injection} (page \pageref{Condition Injection}),
an above injection constructs a Yablo-like structure in $ \xbG,$ by
Definition \ref{Definition Extension} (page \pageref{Definition Extension}),
we have a full valuation of $ \xbG,$ which preserves this structure.
\subsection{
If a Suitable Graph Contains Some $x$ Such That $x$ and $\xCN x$
Are Contradictory, Then it Contains the Yablo Graph.
}

% {\tiny LABEL: {Section Converse}} \\[1mm]
\label{Section Converse}

\efa

\bfa

$\hspace{0.01em}$

% {\tiny (+++ Orig. No.:  Fact Simplest +++)}

% {\tiny LABEL: {Fact Simplest}} \\[1mm]
\label{Fact Simplest}

Yablo's structure is the simplest structure showing that both $x_{0}$ and
$ \xCN x_{0}$
are contradictory, and contained in all such structures.

\efa

\subparagraph{
Proof
}

$\hspace{0.01em}$

% {\tiny (+++ Orig.:  Proof +++)}

Earlier sections of this paper. More precisely:
 \xEh
 \xDH The basic strategy (contradictions, and satisfying antagonistic
requirements) is without alternatives.
See Section \ref{Section Strategy} (page \pageref{Section Strategy}).
 \xDH The elementary contradiction cells have to be (essentially) as in
Yablo's construction. In particular, they have to be triangles,
with all three sides negative.
See Section 
\ref{Section Elementary-Cells} (page 
\pageref{Section Elementary-Cells}).
 \xDH Combining contradiction cells may be somewhat different, but they
have to contain Yablo's structure.
See Section 
\ref{Section Combining-Cells} (page 
\pageref{Section Combining-Cells}),
in particular Section 
\ref{Section Combining-Cells} (page 
\pageref{Section Combining-Cells}),
Remark \ref{Remark Arbitrary} (page \pageref{Remark Arbitrary}).
See also the Saw Blade construction e.g. in
 \cite{Sch22} or  \cite{Sch23b}.
 \xEj

\bfa

$\hspace{0.01em}$

% {\tiny (+++ Orig. No.:  Fact Xi +++)}

% {\tiny LABEL: {Fact Xi}} \\[1mm]
\label{Fact Xi}

Suppose the graph $ \xbG $ has a valuation with an element $ \xfA x_{0}
\xfA = \xbx.$
Then, we can reconstruct Yablo's structure within $ \xbG $ (and thus also
obtain
an injection as described in
Condition \ref{Condition Injection} (page \pageref{Condition Injection})).

\efa

\subparagraph{
Proof
}

$\hspace{0.01em}$

% {\tiny (+++ Orig.:  Proof +++)}

By Fact \ref{Fact Simplest} (page \pageref{Fact Simplest}), or
directly by the construction
Construction \ref{Construction Details} (page \pageref{Construction Details})
inside $ \xbG,$ we have an analogue of the Yablo structure in $ \xbG,$
perhaps with paths instead of simple arrows, and perhaps with
many unnecessary parts.

We can eliminate the unnecessary parts as follows, (for ease of reading
with a simple example, analogously to
Definition \ref{Definition Extension} (page \pageref{Definition Extension})):

Suppose we have $x_{i}=x_{i+1} \xcu x_{i}',$ where the part $x_{i}' $ is
unnecessary.
Then replace $x_{i}' $ by $(x_{i}' \xco \xCN x_{i}' \xco \xdT),$ so we
have $x_{i}=x_{i+1} \xcu (x_{i}' \xco \xCN x_{i}' \xco \xdT).$

Thus, we have the same structure as the one defined by
Condition \ref{Condition Injection} (page \pageref{Condition Injection})  and
Definition \ref{Definition Extension} (page \pageref{Definition Extension}).
\subsection{
Conclusion
}

We summarize:

The existence of a possible substructure of $ \xbG $ according to
Condition \ref{Condition Injection} (page \pageref{Condition Injection})  and
Definition \ref{Definition Extension} (page \pageref{Definition Extension})
is equivalent to the existence of an element $x$ with $ \xfA x \xfA = \xbx
.$
\clearpage
\section{
Comparison to related work
}

% {\tiny LABEL: {Section Comparison}} \\[1mm]
\label{Section Comparison}
\subsection{
Comments on Rabern et al., \cite{RRM13}
}

% {\tiny LABEL: {Section Rabern}} \\[1mm]
\label{Section Rabern}
\subsubsection{
Logics in Rabern et al.
}

% {\tiny LABEL: {Section Logics-Rabern}} \\[1mm]
\label{Section Logics-Rabern}

The ideas in our more informal comments on a logic for reasoning in graphs
in
Section 
\ref{Section Logic} (page 
\pageref{Section Logic})  were more elaborated in
 \cite{RRM13}, Section 1 there, e.g. our notion of ``model''
corresponds
to ``acceptable truth value assignment'' there.

Our discussion in
Section \ref{Section Rules} (page \pageref{Section Rules}), however,
and the allusion to problems with multiple origins,
see Example \ref{Example Logic} (page \pageref{Example Logic}), (7),
seem to be missing in  \cite{RRM13}.
\subsubsection{
A conjecture by Rabern et al.
}

% {\tiny LABEL: {Section Conjecture-Rabern}} \\[1mm]
\label{Section Conjecture-Rabern}

This result was published e.g. in
 \cite{Sch21}, v1, and  \cite{Sch22}.

We show that conjecture 15 in
 \cite{RRM13} is wrong.

\be

$\hspace{0.01em}$

% {\tiny (+++ Orig. No.:  Example YG' +++)}

% {\tiny LABEL: {Example YG'}} \\[1mm]
\label{Example YG'}

We define now a modified Yablo graph $ \xCf YG',$ and a corresponding
denotation $d,$
which is paradoxical.

We refer to Fig.3 in  \cite{RRM13}, and
Diagram \ref{Diagram YG'} (page \pageref{Diagram YG'}).

 \xEh
 \xDH
The vertices (and the set $S$ of language symbols):

We keep all $Y_{i}$ of Fig.3 in  \cite{RRM13},
and introduce new vertices $(Y_{i},Y_{j},Y_{k})$ for $i<k<j.$
(When we write $(Y_{i},Y_{j},Y_{k}),$ we tacitly assume that $i<k<j.)$
 \xDH
The arrows:

All $Y_{i} \xcp Y_{i+1}$ as before.
We ``factorize'' longer arrows through new vertices:
 \xEh
 \xDH
$Y_{i} \xcp (Y_{i},Y_{j},Y_{i+1})$
 \xDH
$(Y_{i},Y_{j},Y_{k}) \xcp (Y_{i},Y_{j},Y_{k+1})$
 \xDH
$(Y_{i},Y_{j},Y_{j-1}) \xcp Y_{j}$
 \xEj
See Diagram \ref{Diagram YG'} (page \pageref{Diagram YG'}).
 \xEj

We define $d$ (instead of writing $d((x,y,z))$ we write $d(x,y,z)$ -
likewise $[x,y,z]_{v}$ for $[(x,y,z)]_{v}$ below):
 \xEh
 \xDH
$d(Y_{i}):= \xCN Y_{i+1} \xcu \xcU \{ \xCN (Y_{i},Y_{j},Y_{i+1}):$ $i+2
\xck j\}$

(This is the main idea of the Yablo construction.)
 \xDH
$d(Y_{i},Y_{j},Y_{k})$ $:=$ $(Y_{i},Y_{j},Y_{k+1})$ for $i<k<j-1$
 \xDH
$d(Y_{i},Y_{j},Y_{j-1})$ $:=$ $Y_{j}$
 \xEj

Obviously, $ \xCf YG' $ corresponds to $S$ and $d,$ i.e. $YG' =G_{S,d}.$

\ee

\bfa

$\hspace{0.01em}$

% {\tiny (+++ Orig. No.:  Fact YG' +++)}

% {\tiny LABEL: {Fact YG'}} \\[1mm]
\label{Fact YG'}

$ \xCf YG' $ and $d$ code the Yablo Paradox:

\efa

\subparagraph{
Proof
}

$\hspace{0.01em}$

% {\tiny (+++ Orig.:  Proof +++)}

Let $v$ be an acceptable valuation relative to $d.$

Suppose $[Y_{1}]_{v}= \xdT,$ then $[Y_{2}]_{v}= \xdF,$ and
$[Y_{1},Y_{k},Y_{2}]_{v}= \xdF $ for $2<k,$ so
$[Y_{k}]_{v}= \xdF $ for $2<k.$
By $[Y_{2}]_{v}= \xdF,$ there must be $j$ such that
$j=3$ and $[Y_{3}]_{v}= \xdT,$ or $j>3$ and $[Y_{2},Y_{j},Y_{3}]_{v}=
\xdT,$ and
again, $[Y_{j}]_{v}= \xdT,$
a contradiction.

If $[Y_{1}]_{v}= \xdF,$ then as above for $[Y_{2}]_{v},$ we find $j \xcg
2$ and $[Y_{j}]_{v}= \xdT,$ and
argue with $Y_{j}$ as above for $Y_{1}.$

Thus, $ \xCf YG' $ with $d$ as above is paradoxical, and $ \xCf YG' $ is
dangerous.

$ \xcz $
\\[3ex]

\be

$\hspace{0.01em}$

% {\tiny (+++ Orig. No.:  Example Homomorphismus +++)}

% {\tiny LABEL: {Example Homomorphismus}} \\[1mm]
\label{Example Homomorphismus}

We first define $ \xCf YG'' $:
$V(YG''):=\{ \xBc Y_{i} \xBe:i< \xbo \},$ $E(YG''):=\{ \xBc Y_{i} \xBe  \xcp 
\xBc Y_{i+1} \xBe:i<
\xbo \}.$

We now define the homomorphism from $ \xCf YG' $ to $ \xCf YG''.$
We collaps for fixed $k$ $Y_{k}$ and all $(Y_{i},Y_{j},Y_{k})$ to
$ \xBc Y_{k} \xBe,$ more precisely,
define $f$
by $f(Y_{k}):=f(Y_{i},Y_{j},Y_{k}):= \xBc Y_{k} \xBe $ for all suitable $i,j.$

Note that $ \xCf YG' $ only had arrows between ``successor levels'',
and we have now only arrows from $ \xBc Y_{k} \xBe $ to $ \xBc Y_{k+1} \xBe,$
so $f$ is a
homomorphism,
moreover, our structure $ \xCf YG'' $ has the
is not dangerous, contradicting
conjecture 15 in  \cite{RRM13}.

\clearpage

\begin{diagram}

% {\tiny LABEL: {Diagram YG'}} \\[1mm]
\label{Diagram YG'}
\index{Diagram YG'}

\unitlength0.65mm
\begin{picture}(150,180)(0,0)

\put(0,175){{\rm\bf Diagram YG' }}

\put(0,160){$Y_1$}
\put(140,160){$ \xBc Y_1 \xBe $}
\put(2,158){\line(0,-1){32}}
\put(143,158){\line(0,-1){33}}
\put(3,158){\line(3,-4){25}}
\put(4,158){\line(4,-3){45}}
\put(5,158){\line(5,-2){80}}

\put(0,120){$Y_2$}
\put(2,118){\line(0,-1){32}}
\put(143,118){\line(0,-1){33}}
\put(20,120){{\footnotesize $(Y_1,Y_3,Y_2)$}}
\put(50,120){{\footnotesize $(Y_1,Y_4,Y_2)$}}
\put(80,120){{\footnotesize $(Y_1,Y_5,Y_2)$}}
\put(140,120){$ \xBc Y_2 \xBe $}
\put(3,118){\line(5,-2){80}}
\put(4,118){\line(3,-1){105}}
\put(27,118){\line(-3,-4){25}}
\put(55,118){\line(-3,-4){25}}
\put(85,118){\line(-3,-4){25}}

\put(0,80){$Y_3$}
\put(2,78){\line(0,-1){32}}
\put(143,78){\line(0,-1){33}}
\put(20,80){{\footnotesize $(Y_1,Y_4,Y_3)$}}
\put(50,80){{\footnotesize $(Y_1,Y_5,Y_3)$}}
\put(80,80){{\footnotesize $(Y_2,Y_4,Y_3)$}}
\put(110,80){{\footnotesize $(Y_2,Y_5,Y_3)$}}
\put(140,80){$ \xBc Y_3 \xBe $}
\put(3,78){\line(5,-2){80}}
\put(27,78){\line(-3,-4){25}}
\put(55,78){\line(-3,-4){25}}
% \put(85,78){\line(-5,-2){80}}
% \put(85,78){\line(-3,-1){80}}
\put(80,79){\line(-2,-1){73}}
\put(115,78){\line(-3,-2){50}}

\put(0,40){$Y_4$}
\put(2,38){\line(0,-1){32}}
\put(143,38){\line(0,-1){33}}
\put(20,40){{\footnotesize $(Y_1,Y_5,Y_4)$}}
\put(50,40){{\footnotesize $(Y_2,Y_5,Y_4)$}}
\put(80,40){{\footnotesize $(Y_3,Y_5,Y_4)$}}
\put(140,40){$ \xBc Y_4 \xBe $}
\put(27,38){\line(-3,-4){25}}
% \put(55,38){\line(-5,-4){45}}
\put(56,38){\line(-3,-2){51}}
% \put(85,38){\line(-5,-2){80}}
% \put(85,38){\line(-3,-1){80}}
\put(80,39){\line(-2,-1){73}}

\put(0,0){$Y_5$}
\put(140,0){$ \xBc Y_5 \xBe $}

% \put(43,27){\vector(1,1){24}}

% \put(53,67){\line(1,1){4}}

\end{picture}

\end{diagram}

\vspace{4mm}

\ee

This is just the start of the graph, it continues downward through
$ \xbo $ many levels.

The lines stand for downward pointing arrows.
The lines originating from the $Y_{i}$ correspond to the negative
lines in the original Yablo graph, all others are simple positive
lines, of the type $d(X)=X'.$

The left part of the drawing represents the graph YG', the
right hand part the collapsed graph, the homomorphic image YG".

Compare to Fig.3 in  \cite{RRM13}.
\clearpage
\subsection{
Comments on Beringer and Schindler, \cite{BS17}
}

% {\tiny LABEL: {Section Beringer}} \\[1mm]
\label{Section Beringer}
\subsubsection{
Introduction
}

We prove a modified version of conjecture 2.4
of  \cite{BS17}.

The proof is mainly based on the minimality of the Yablo structure,
the additions are often reformulations between different situations.

\paragraph{
Conjecture 2.4 in  \cite{BS17}
}

$\hspace{0.01em}$

% {\tiny (+++ Orig.:  Conjecture 2.4 in  [BS17 ] +++)}

% {\tiny LABEL: {Section Conjecture 2.4 in  [BS17 ]}} \\[1mm]
\label{Section Conjecture 2.4 in  [BS17 ]}

(Conjecture 2.4 in  \cite{BS17}, $p.$ 474):
A reference graph is dangerous iff it contains a subdivision of the
liar-graph as a subgraph or the Yablo-graph as a finitary minor.

As we consider only directed, acyclic graphs, and whether they permit a
Yablo-like interpretation, we modify to:

 \xEI
 \xDH An acyclic graph has an element $x$ such that assuming $x$, as well
as assuming
$ \xCN x$ lead to a contradiction (we call this the Yablo property)
iff it contains the Yablo-graph as a finitary minor.
 \xEJ

We repeat here for the reader's convenience

\bd

$\hspace{0.01em}$

% {\tiny (+++ Orig. No.:  Definition 6.4 +++)}

% {\tiny LABEL: {Definition 6.4}} \\[1mm]
\label{Definition 6.4}

(Definition 6.4 in  \cite{BS17}, $p.$ 490):

Call a directed graph $G$ a finitary inflation of a directed graph $H$
if $V(G)$ admits a partition $\{V_{x}:x \xbe V(H)\}$ such that the
following
conditions hold:

1. for all $x,y \xbe V(H):$

(a) if $x \xEd y,$ then $(x,y) \xbe A(H)$ iff there is some $a \xbe A(G)$
with tail in $V_{x}$ and
head in $V_{y},$

(b) if $x=y,$ then $(x,y) \xce A(H),$

2. for all $x \xbe V(H)$ and all $y,z \xbe V_{x} \xcc V(G):$ if $y$ has an
in-neighbour in $V(G)-V_{x}$
and $z$ has an out-neighbour in $V(G)-V_{x}$
then there is a path $P \xcc G[V_{x}]$ from $y$ to $z,$

3. for all $x \xbe V(H),$ $V_{x}$ is finite.

$(V(G),$ $A(G),$ tail, head, in-neighbour, out-neighbour etc. are the
obvious.)

\ed

\bd

$\hspace{0.01em}$

% {\tiny (+++ Orig. No.:  Definition 6.5 +++)}

% {\tiny LABEL: {Definition 6.5}} \\[1mm]
\label{Definition 6.5}

(Definition 6.5 in  \cite{BS17}, $p.$ 490):

A directed graph $H$ is a finitary minor of $G' $ iff there is some $G
\xcc G' $ such that
$G$ is a finitary inflation of $H.$

\clearpage

\begin{diagram}

% {\tiny LABEL: {Diagram BS}} \\[1mm]
\label{Diagram BS}
\index{Diagram BS}

\unitlength0.8mm
\begin{picture}(220,200)(0,0)

\put(40,40){\circle*{1}}
\put(100,40){\circle*{1}}
\put(160,40){\circle*{1}}

\put(40,40){\line(1,2){25}}
\put(40,40){\line(-1,2){25}}
\put(100,40){\line(1,2){25}}
\put(100,40){\line(-1,2){25}}
\put(160,40){\line(1,2){25}}
\put(160,40){\line(-1,2){25}}

\put(38,35){\footnotesize $V_1$}
\put(98,35){\footnotesize $V_2$}
\put(158,35){\footnotesize $V_3$}

\put(10,70){\line(1,0){18}}
\put(51,70){\line(1,0){37}}
\put(51,80){\line(1,0){37}}
\put(111,70){\line(1,0){38}}
\put(172,70){\line(1,0){18}}

\put(12,65){\footnotesize $a_{01}$}
\put(68,65){\footnotesize $a_{12}$}
\put(68,77){\footnotesize $a''_{12}$}
\put(128,65){\footnotesize $a_{23}$}
\put(182,65){\footnotesize $a_{34}$}

\put(200,70){$G$}
\put(200,22){$H$}

\put(50,100){\line(0,-1){10}}
\put(50,100){\line(1,0){100}}
\put(150,100){\line(0,-1){10}}
\put(102,96){\footnotesize $a_{13}$}

\multiput(30,70)(1,0){20}{\circle*{0.3}}
\multiput(30,70)(2,1){10}{\circle*{0.3}}
\multiput(30,70)(1,1){20}{\circle*{0.3}}
\multiput(90,70)(1,0){20}{\circle*{0.3}}
\multiput(110,70)(-2,1){10}{\circle*{0.3}}
\multiput(170,70)(-1,0){20}{\circle*{0.3}}
\multiput(170,70)(-1,1){20}{\circle*{0.3}}

\put(38,66){\footnotesize $ \xbs _{1}$}
\put(38,86){\footnotesize $ \xbs' _{1}$}
\put(45,74){\footnotesize $ \xbs'' _{1}$}
\put(98,66){\footnotesize $ \xbs _{2}$}
\put(99,77){\footnotesize $ \xbs'' _{2}$}
\put(158,66){\footnotesize $ \xbs _{3}$}
\put(158,86){\footnotesize $ \xbs' _{3}$}

\put(40,25){\circle*{1}}
\put(100,25){\circle*{1}}
\put(160,25){\circle*{1}}

\put(42,25){\line(1,0){56}}
\put(102,25){\line(1,0){56}}
\put(40,10){\line(0,1){9}}
\put(40,10){\line(1,0){120}}
\put(160,10){\line(0,1){9}}

\put(38,20){\footnotesize $x_1$}
\put(98,20){\footnotesize $x_2$}
\put(158,20){\footnotesize $x_3$}

\end{picture}

\end{diagram}

\vspace{4mm}

\ed

\br

$\hspace{0.01em}$

% {\tiny (+++ Orig. No.:  Remark Even +++)}

% {\tiny LABEL: {Remark Even}} \\[1mm]
\label{Remark Even}

This (trivial) remark underlines the importance of free choice of
values for arrows.

Consider $x_{0},$ $x_{1},$ $x_{2},$ $x_{3},$ $x_{4},$ (as in Yablo's
construction), and
$V_{x_{0}},$  \Xl., $V_{x_{4}},$ abbreviated $V_{0},$ etc.

In Yablo's construction, we have $x_{1} \xcP x_{2} \xcP x_{3},$ $x_{1}
\xcP x_{3}$ (negated arrows).

By above Definition, 1. (a),
we have arrows $a_{12}:V_{1} \xcp V_{2},$ and perhaps some additional
$a_{12}'':V_{1} \xcp V_{2},$ and
$a_{23}:V_{2} \xcp V_{3},$ $a_{13}:V_{1} \xcp V_{3},$ as well as arrows
$a_{01}:V_{0} \xcp V_{1},$ $a_{34}:V_{3} \xcp V_{4}.$

Let, for arrows $a_{ij}$ and paths $ \xbs_{ij}$ be $s(a_{ij})$ the start
of $a_{ij},$ $e(a_{ij})$
its end, likewise for $ \xbs_{ij}.$

By above Definition, 2.,
we have paths $ \xbs_{1}:e(a_{01}) \xcp s(a_{12}),$ $ \xbs_{1}''
:e(a_{01}) \xcp s(a_{12}''),$
$ \xbs'_{1}:e(a_{01}) \xcp s(a_{13}),$
$ \xbs_{2}:e(a_{12}) \xcp s(a_{23}),$ $ \xbs''_{2}:e(a_{12}'') \xcp
s(a_{23}),$
$ \xbs_{3}:e(a_{23}) \xcp s(a_{34}),$ $ \xbs'_{3}:e(a_{13}) \xcp
s(a_{34}).$

Suppose $length(\xbs_{2})$ is odd, the length of all other paths $
\xbs_{ij}$ etc. is
even.

Then there is no contradiction, if all arrows have to be negative.

See Diagram \ref{Diagram BS} (page \pageref{Diagram BS}).
(The diagonal lines mark the limits of the $V_{i},$ Other solid lines
stand
for arrows from left to right, dotted lines stand for paths from left
to right.)

\er

Before we begin with the proofs, we postulate a simplification:

\bd

$\hspace{0.01em}$

% {\tiny (+++ Orig. No.:  Definition Canonical Partition +++)}

% {\tiny LABEL: {Definition Canonical Partition}} \\[1mm]
\label{Definition Canonical Partition}

(Canonical partition)

In Yablo's structure, every node has finitely many predecessors, but
infinitely successors. For this reason, and finiteness of all $V_{x_{j}},$
we may assume that every element
of the partition $V_{x_{j}}$ has as many ``entry points'' $y_{i,j}$ as
predecessors $V_{x_{i}},$
but only one exit point $z_{j}.$ $V_{x_{i}}$ and $V_{x_{j}}$ are joined by
one arrow $z_{i} \xcP y_{i,j},$
and there is a path $ \xbs_{i,j}:y_{i,j} \Xl z_{j}.$

See Diagram \ref{Diagram Partition} (page \pageref{Diagram Partition}).
\subsubsection{
Proof of the conjecture from partition to Yablo property.
}

\ed

We give the start of the Yablo structure as example, the generalization
is obvious.

\be

$\hspace{0.01em}$

% {\tiny (+++ Orig. No.:  Example BS-Structure +++)}

% {\tiny LABEL: {Example BS-Structure}} \\[1mm]
\label{Example BS-Structure}

Consider as in Yablo $x_{0} \xcP x_{1} \xcP x_{2} \xcP x_{3},$ $x_{0} \xcP
x_{2},$ $x_{0} \xcP x_{3},$ $x_{1} \xcP x_{3}$ with arrows

$b_{0,1}:x_{0} \xcP x_{1},$ $b_{0,2}:x_{0} \xcP x_{2},$ $b_{0,3}:x_{0}
\xcP x_{3},$
$b_{1,2}:x_{1} \xcP x_{2},$ $b_{1,3}:x_{1} \xcP x_{3},$ $b_{2,3}:x_{2}
\xcP x_{3},$
etc. in Yablo Structure.

(The $H$ in Diagram \ref{Diagram BS} (page \pageref{Diagram BS}).)

We then have as $G$ (notation as in
Diagram \ref{Diagram Partition} (page \pageref{Diagram Partition}))
the partition:

 \xEh
 \xDH
$V_{x_{0}},$ $V_{x_{1}},$ $V_{x_{2}},$ $V_{x_{3}},$ and

 \xDH (negative) arrows

$a_{0,1}:z_{0} \xcP y_{0,1},$ $a_{0,2}:z_{0} \xcP y_{0,2},$ $a_{0,3}:z_{0}
\xcP y_{0,3},$
$a_{1,2}:z_{1} \xcP y_{1,2},$ $a_{1,3}:z_{0} \xcP y_{1,3},$ $a_{2,3}:z_{0}
\xcP y_{2,3},$

 \xDH (positive) paths

$ \xbs_{0,1}:y_{0,1} \Xl z_{1},$ $ \xbs_{0,2}:y_{0,2} \Xl z_{2},$ $
\xbs_{0,3}:y_{0,3} \Xl z_{3},$
$ \xbs_{1,2}:y_{1,2} \Xl z_{2},$ $ \xbs_{1,3}:y_{1,3} \Xl z_{3},$ $
\xbs_{2,3}:y_{2,3} \Xl z_{3},$

 \xEj

Obviously, the structure has the Yablo property, the paths $ \xbs_{i,j}$
work as neutral elements.
\subsubsection{
Construction of a partition from a graph with the Yablo property.
}

\ee

First, a restriction on the graphs we consider:

\bd

$\hspace{0.01em}$

% {\tiny (+++ Orig. No.:  Definition Additional Prerequisite +++)}

% {\tiny LABEL: {Definition Additional Prerequisite}} \\[1mm]
\label{Definition Additional Prerequisite}

(Additional Prerequisite)

We work with paths here, but we do not permit branchings as follows:
$ \xCf x \Xl y \Xl z,$ with a branch continuing $ \xCf y \Xl z'.$ We
impose branchings at the
beginning, i.e. $ \xCf x \Xl y \Xl z$, and $x \Xl y'  \Xl z',$ they run
parallel until $y,$ and
then diverge. We need this, as we will make the first arrow decide where
we go. This seems a minor modification.

The reason for this is: we may have recursive branchings of arbitrarily
finite size, and it does not seem obvious how to find the ``right''
arrow. The best way to circumvent this problem is to make some kind
of preprocessing. (This is implicitely packed in our approach into the
Condition \ref{Condition Injection} (page \pageref{Condition Injection}).)

\ed

Assume there is a graph with the Yablo Property YP, by minimality, it
contains
the Yablo structure, perhaps with paths.
(We assume that the additional prerequisite holds.)

Again, we give an example, the generalization is again obvious.

\be

$\hspace{0.01em}$

% {\tiny (+++ Orig. No.:  Example Yablo-Property +++)}

% {\tiny LABEL: {Example Yablo-Property}} \\[1mm]
\label{Example Yablo-Property}

We have a graph with the Yablo property, perhaps formed by paths instead
of single arrows, but with the additional prerequisites.

By minimality, we have some

$ \xbp_{0,1}:x_{0} \xcP x_{1},$ $ \xbp_{0,2}:x_{0} \xcP x_{2},$ $
\xbp_{0,3}:x_{0} \xcP x_{3},$
$ \xbp_{1,2}:x_{1} \xcP x_{2},$ $ \xbp_{1,3}:x_{1} \xcP x_{3},$ $
\xbp_{2,3}:x_{2} \xcP x_{3},$

which decompose into

$ \xbp_{0,1}=a_{0,1} \xDO \xbs_{0,1},$ $ \xbp_{0,2}=a_{0,2} \xDO
\xbs_{0,2},$ $ \xbp_{0,3}=a_{0,3} \xDO \xbs_{0,3},$
$ \xbp_{1,2}=a_{1,2} \xDO \xbs_{1,2},$ $ \xbp_{1,2}=a_{1,3} \xDO
\xbs_{1,3},$ $ \xbp_{2,3}=a_{2,3} \xDO \xbs_{2,3}.$

(Note that they diverge on the first arrows.)

We construct a partition as follows:

 \xEh
 \xDH
$V_{x_{0}},$ $V_{x_{1}},$ $V_{x_{2}},$ $V_{x_{3}},$ and

 \xDH (negative) arrows

$a_{0,1}:z_{0} \xcP y_{0,1},$ $a_{0,2}:z_{0} \xcP y_{0,2},$ $a_{0,3}:z_{0}
\xcP y_{0,3},$
$a_{1,2}:z_{1} \xcP y_{1,2},$ $a_{1,3}:z_{0} \xcP y_{1,3},$ $a_{2,3}:z_{0}
\xcP y_{2,3},$

 \xDH (positive) paths

$ \xbs_{0,1}:y_{0,1} \Xl z_{1},$ $ \xbs_{0,2}:y_{0,2} \Xl z_{2},$ $
\xbs_{0,3}:y_{0,3} \Xl z_{3},$
$ \xbs_{1,2}:y_{1,2} \Xl z_{2},$ $ \xbs_{1,3}:y_{1,3} \Xl z_{3},$ $
\xbs_{2,3}:y_{2,3} \Xl z_{3},$

 \xEj

\clearpage

\begin{diagram}

% {\tiny LABEL: {Diagram Partition}} \\[1mm]
\label{Diagram Partition}
\index{Diagram Partition}

\unitlength0.8mm
\begin{picture}(250,300)(0,0)

\put(40,40){\circle*{1}}
\put(120,40){\circle*{1}}
\put(200,40){\circle*{1}}

\put(38,35){\footnotesize $V_{x_i}$}
\put(118,35){\footnotesize $V_{x_j}$}
\put(198,35){\footnotesize $V_{x_k}$}

\thicklines

\put(40,40){\line(1,8){30}}
\put(40,40){\line(-1,8){30}}
\put(120,40){\line(1,8){30}}
\put(120,40){\line(-1,8){30}}
\put(200,40){\line(1,8){30}}
\put(200,40){\line(-1,8){30}}

\put(40.5,40){\line(1,8){30}}
\put(39.5,40){\line(-1,8){30}}
\put(120.5,40){\line(1,8){30}}
\put(119.5,40){\line(-1,8){30}}
\put(200.5,40){\line(1,8){30}}
\put(199.5,40){\line(-1,8){30}}

\put(50,120){\line(1,0){140}}
\put(50,120){\line(21,24){52.5}}
\put(137.5,180){\line(1,0){45}}

\thinlines

\put(50,120){\circle*{1}}
\put(42,120){\footnotesize $z_{i}$}

\put(115,80){\circle*{1}}
\put(105,80){\footnotesize $y_{0,j}$}
\put(102.5,180){\circle*{1}}
\put(95,185){\footnotesize $y_{i,j}$}
\put(137.5,180){\circle*{1}}
\put(132,185){\footnotesize $z_{j}$}

\put(195,80){\circle*{1}}
\put(185,80){\footnotesize $y_{0,k}$}
\put(190,120){\circle*{1}}
\put(183,115){\footnotesize $y_{i,k}$}
\put(182.5,180){\circle*{1}}
\put(176,173){\footnotesize $y_{j,k}$}
\put(222.5,220){\circle*{1}}
\put(225,220){\footnotesize $z_{k}$}

\put(102.5,180){\line(1,0){35}}

\put(70,150){\footnotesize $a_{i,j}$}
\put(160,115){\footnotesize $a_{i,k}$}
\put(160,185){\footnotesize $a_{j,k}$}

\put(115,80){\line(9,40){22.5}}

\put(182.5,180){\line(1,1){40}}

\put(190,120){\line(13,40){32.5}}

\put(195,80){\line(11,56){27.5}}

\put(115,183){\footnotesize $ \xbs _{i,j}$}
\put(118,130){\footnotesize $ \xbs _{0,j}$}

\put(200,100){\footnotesize $ \xbs _{0,k}$}
\put(200,170){\footnotesize $ \xbs _{i,k}$}
\put(195,200){\footnotesize $ \xbs _{j,k}$}

% \put(38,35){\footnotesize $V_1$}

% \put(12,65){\footnotesize $a_{01}$}

% \put(200,70){$G$}
% \put(200,22){$H$}

% \multiput(30,70)(1,0){20}{\circle*{0.3}}

% \put(38,66){\footnotesize $ \xbs _{1}$}

\end{picture}

\end{diagram}

\vspace{4mm}

\clearpage

\ee

.
\clearpage

$ \xCO $
% {\tiny KS= wab-Main   FOUND IN END
% c:-ks-w-t-framart-phil-yablo-wablo-wbe-1.mp} \\[1mm]

\end{document}